\documentclass{amsart}

\newtheorem{theorem}{Theorem}[section]
\newtheorem{lemma}[theorem]{Lemma}
\newtheorem{corollary}[theorem]{Corollary}     
\newtheorem*{maintheorem}{Main Theorem}
\newtheorem*{maintheorem_alt}{Main Theorem (Alternative Form)}

\theoremstyle{definition}
\newtheorem{definition}[theorem]{Definition}

\theoremstyle{remark}

\numberwithin{equation}{section}

\newcommand{\boldC}{{\bf C}}
\newcommand{\boldD}{{\bf D}}
\newcommand{\boldN}{{\bf N}}
\newcommand{\boldS}{{\bf S}}

\newcommand{\eE}{{\mathcal E}}      
\newcommand{\eG}{{\mathfrak G}}
\newcommand{\eH}{{\mathfrak H}}
\newcommand{\eK}{{\mathfrak K}}
\newcommand{\eL}{{\mathfrak L}}
\newcommand{\eM}{{\mathfrak M}}
\newcommand{\eN}{{\mathfrak N}}
\newcommand{\eX}{{\mathfrak X}}

\newcommand{\calD}{{\mathcal D}}
\newcommand{\calE}{{\mathcal E}}
\newcommand{\calF}{{\mathcal F}}
\newcommand{\calV}{{\mathcal V}}

\newcommand{\tildeC}{\widetilde C}
\newcommand{\Skappa}{\boldS_\kappa}
\newcommand{\SkappaCplus}{\boldS_\kappa(\boldC_+)}
\newcommand{\Szero}{\boldS_0}
\newcommand{\Nkappa}{\boldN_\kappa}
\newcommand{\Nzero}{\boldN_0}
\newcommand{\Skappaprime}{\boldS_{\kappa'}}
\newcommand{\Skappadoubleprime}{\boldS_{\kappa''}}
\newcommand{\Nkappaprime}{\boldN_{\kappa'}}
\newcommand{\ip}[2]{\left<#1,#2\right>}
\newcommand{\indminus}{\hbox{\rm{ind}}_-\,}
\newcommand{\sqminus}{{\rm{sq}_-}}
\newcommand{\oh}{{\mathcal O}}
\newcommand{\lh}{\eL(\eH)}
\newcommand{\diag}{\hbox{\rm{diag}}\,}

\begin{document}

\title[Notes on Interpolation:  Nudel$'$man's Problem]
{Notes on Interpolation in the Generalized  \\ [8pt]
   Schur Class. II.  Nudel$'$man's Problem} 

\author{D. Alpay}
\author{T. Constantinescu} 
\author{A. Dijksma} 
\author{J. Rovnyak}
\address{Department of Mathematics \\
  Ben-Gurion University of the Negev \\
  P.O. Box 653 \\
  84105 Beer-Sheva, Israel} 
\email{\tt dany@ivory.bgu.ac.il}
\address{Programs in Mathematical Sciences \\
  University of Texas at Dallas \\
  Box 830688, Richardson, TX 75083-0688, U. S. A.}
\email{tiberiu@utdallas.edu} 
\address{
  Department of Mathematics \\
  University of Groningen \\
  P.O. Box 800 \\
  9700 AV Groningen, The Netherlands} 
\email{\tt dijksma@math.rug.nl}
\address{\!\!\! University of Virginia,
  \!\! Department of Mathematics \\
  P.O. \!Box 400137 \\
  \!Charlottesville, VA 22904-4137, U. S. A.}  
\email{\tt
  rovnyak@Virginia.EDU}
\date{} 
\thanks{
  {\it $2000$ Mathematics Subject Classifications}.  
    Primary 47A57 30E05 47B32.  Secondary 47B50 42A50 \\
  \indent \hskip.06cm   
   J. Rovnyak was supported by
  the National Science Foundation  DMS-0100437 and by the
  Netherlands Organization for Scientific Research NWO~B~61-482.}

\begin{abstract}
  An indefinite generalization of Nudel$'$man's problem is used in a
  systematic approach to interpolation theorems for generalized Schur
  and Nevanlinna functions with interior and boundary data.  Besides
  results on existence criteria for Pick-Nevanlinna and
  Carath\'eodory-Fej\'er interpolation, the method yields new results
  on generalized interpolation in the sense of Sarason and boundary
  interpolation, including properties of the finite Hilbert transform
  relative to weights.  The main theorem appeals to the Ball and
  Helton almost-commutant lifting theorem to provide criteria for the
  existence of a solution to Nudel$'$man's problem.
\end{abstract}

\maketitle \setcounter{equation}{0}

\section{Introduction}


In this paper we use a theorem of Ball and Helton \cite{BH} to
describe a unified approach to a series of classical interpolation
problems for the generalized Schur and Nevanlinna classes.  These
classes contain meromorphic functions having a finite number of poles
on a disk or half-plane.  The main results of the paper update
Chapter~2 of \cite{RR1985}, where similar results are derived for
holomorphic functions.

In \cite{RR1985}, classical interpolation problems are viewed as
special cases of Nudel$'$man's problem: given vectors $b$ and $c$ in a
complex vector space $\calV$ and an operator $A$ on $\calV$ it is
required to find a holomorphic function $f(z)$ which is bounded by one
on the unit disk such that
\begin{equation*}
  b = f(A)c
\end{equation*}
with a suitable interpretation of the equation.  Particular choices of
$A,b$, and $c$ give various classical problems, such as the
Pick-Nevanlinna and Carath\'eodory-Fej\'er problems and boundary
problems of Loewner type.  Conditions for the existence of a solution
of Nudel$'$man's problem are derived in \cite[Chapter 2]{RR1985} as an
application of the theory of contraction operators on a Hilbert space.
This yields existence criteria for classical interpolation problems
for Schur and Nevanlinna functions, that is, holomorphic functions
which are bounded by one on the unit disk or which have nonnegative
imaginary part on the upper half-plane.

Our main tool is an indefinite generalization of Nudel$'$man's
problem, which provides a similar framework for interpolation problems
involving the classes $\Skappa$ and $\Nkappa$ of generalized Schur and
Nevanlinna functions on the unit disk and upper half-plane.  These
well-known classes consist of meromorphic functions having a finite
number of poles.  As in the definite case, we derive existence
criteria for the solvability of the generalized Nudel$'$man problem in
an abstract setting, and the applications to classical problems follow
as specializations of this result.

Classical sources for interpolation in the generalized Schur and
Nevanlinna classes may be found in \cite{AAK}.  The list below
identifies some more recent works that treat these topics and related
areas.  The list is not a complete bibliography, and the authors
regret omissions.  \smallskip

\noindent{Interpolation of interior data}: \cite{AAK, ACDR1, 
bolotnikovdym, dlnp,
FFK, golinskii, hsw, ivanchenko, aanudelman, alsakhnovich}

\noindent{Loewner theory and boundary data}:
 \cite{adl_np2, ar, BH}

\noindent{Parametrization of solutions}:
\cite{AAK, BH, golinskii}

\noindent{Interpolation in the Stieltjes class, including
parametrizations}:
\cite{abdstieltjes, derkach}

\noindent{Moment problems}:
\cite{CG1995, CG1998, KL77}
\smallskip


Our generalization of Nudel$'$man's problem and Main Theorem are
formulated in Section~\ref{S:maintheorem}.
Sections~\ref{S:applications} and~\ref{S:halfplane} contain the
applications of the Main Theorem to classical interpolation problems
in the disk and half-plane cases, respectively; many of these results
parallel the definite case \cite{RR1985}.  The proof of the Main
Theorem is given in Section~\ref{S:proof}.  This paper can be read
independently of \cite{ACDR1}.

\section{Main Theorem}\label{S:maintheorem}

By a {\bf kernel} or {\bf Hermitian form} we mean a complex-valued
function $K(\zeta,z) = \overline{K(z,\zeta)}$ on a product set
$\Omega\times\Omega$.  Such a kernel has $\kappa$ {\bf negative
squares}, in symbols $\sqminus K = \kappa$, if every selfadjoint
matrix $\left( K(\zeta_j,\zeta_i) \right)_{i,j=1}^n$,
$\zeta_1,\dots,\zeta_n\in \Omega$, $n = 1,2,3,\dots$, has at most
$\kappa$ negative eigenvalues, and one such matrix has exactly
$\kappa$ negative eigenvalues (counting multiplicity).  When this
condition is satisfied with $\kappa = 0$, the kernel is said to be
{\bf nonnegative}.  A familiar example of a Hermitian form is the
inner product ${\ip{Tf}{g}}_{\eH}$, where $T$ is a bounded selfadjoint
operator on a Hilbert space $\eH$ and $f$ and $g$ are arbitrary
vectors in~$\eH$.  We say that $T$ is {\bf nonnegative} or has {\bf
$\kappa$ negative squares} according as the associated Hermitian form
has the same property.

The {\bf Pick class} or {\bf Nevanlinna class} is the set $\mathcal P$
of holomorphic functions which have nonnegative imaginary part on the
open upper half-plane $\boldC_+$.  The {\bf Schur class} is the set
$\mathcal S$ of holomorphic functions which are bounded by one on the
open unit disk~$\boldD$.  Let $\boldC$ be the complex plane.
Following Kre{\u\i}n and Langer \cite{KL77}, for any nonnegative
integer~$\kappa$ we define the {\bf generalized Nevanlinna class} as
the set $\Nkappa$ of functions $f(z)$ which are holomorphic on some
subregion $\Omega$ of $\boldC_+$ such that the kernel
$[f(z)-\overline{f(\zeta)}]/(z-\bar \zeta)$ has $\kappa$ negative
squares on $\Omega\times\Omega$.  The {\bf generalized Schur class} is
the set $\Skappa$ of functions $S(z)$ which are holomorphic on some
subregion $\Omega$ of $\boldD$ such that the kernel
$[1-S(z)\overline{S(\zeta)}]/(1-z\bar \zeta)$ has $\kappa$ negative
squares on $\Omega\times\Omega$.  Functions in the generalized
Nevanlinna and Schur classes have analytic continuations to $\boldC_+$
and $\boldD$, respectively, excluding at most $\kappa$ poles, and the
corresponding kernels for the extensions also have $\kappa$ negative
squares.  We understand that such functions are identified with their
meromorphic continuations.  With this convention, $\Nzero = {\mathcal
  P}$ and $\Szero = {\mathcal S}$.  By a theorem of Kre{\u\i}n and
Langer \cite[Theorem 3.2 on p.~382]{KL72}, $\Skappa$ coincides with
the set of functions
\begin{equation}\label{eq:kreinlanger}
S(z) = f(z)/B(z),\qquad B(z) = 
c \prod_{j=1}^\kappa\; \frac{z-a_j}{1 - z\bar a_j},
\end{equation}
where $f(z)$ is in $\Szero$, $a_1,\dots,a_\kappa$ are (not
necessarily distinct) points in $\boldD$, $|c|=1$, and $f(a_j)\neq
0$ for all $j=1,\dots,\kappa$.  A function $B(z)$ on the unit disk of
the preceding form is a {\bf Blaschke product of degree~$\kappa$}.  We
call \eqref{eq:kreinlanger} a {\bf Kre{\u\i}n-Langer factorization}
of a function $S \in \Skappa$.  The factors $f(z)$ and $B(z)$ in a
Kre{\u\i}n-Langer factorization are determined up to constants of
modulus one.

The original inspiration for what is called Nudel$'$man's problem in
\cite[p.~22]{RR1985} is Nudel$'$man \cite{nudelmanproblem}.  We now
consider the following variant.

\noindent {\bf Nudel$'$man's Problem for $\Skappa$}.  
{\it Given vectors $b,c$ in a complex vector space $\calV$ and a
  linear operator $A$ on $\calV$ into itself, find a pair $(f,B)$,
  where $f \in \Szero$ and $B$ is a Blaschke product of
  degree~$\kappa$, such that
  \begin{equation}
    \label{eq:nud2}
    f(A)c = B(A)b.
  \end{equation}
We call $(A,b,c)$ the {\bf data} of the problem.}

Here and in the rest of this section we assume that we are given a
complex vector space $\calV$ with algebraic dual $\calV'$ (no
topologies on these spaces are required).  If $x'\in \calV'$, we write
$(x,x')$ for $x'(x)$.  If $A:\calV \to \calV$ is a linear operator,
its dual $A'$ is the operator $A':\calV' \to \calV'$ such that
$(x,A'x') = (Ax,x')$ for all $x\in\calV$ and $x'\in\calV'$.

\begin{definition}\label{D:admissible}
  A set $\calD \subseteq \calV'$ is {\bf admissible} for given data
  $(A,b,c)$ if
  \begin{enumerate}
  \item[(i)] $\calD$ is a linear subspace of $\calV'$ which is
    invariant under $A'$, \vskip2pt
  \item[(ii)] the sums $\sum_{j=0}^\infty |(A^jb,x')|^2$ and
    $\sum_{j=0}^\infty |(A^jc,x')|^2$ are finite for all $x' \in
    \calD$, and
  \item[(iii)] if some $x'$ in $\calD$ annihilates $c,Ac,A^2c,\dots$,
    then $x'$ annihilates $b,Ab,A^2b,\dots\;$.
  \end{enumerate}
\end{definition}

Throughout we write $H^2$ for the Hardy space on the unit disk.

\begin{definition}\label{nudelmaneq}
  If $\calD$ is admissible for $(A,b,c)$, then for any
  \smash{$f(z)= \sum_{j=0}^\infty f_j z^j$} and \smash{$B(z)=
    \sum_{j=0}^\infty B_j z^j$} in $H^2$, we interpret \eqref{eq:nud2}
  to mean that
\begin{equation}\label{eq:twosums}
  \sum_{j=0}^\infty f_j (A^j c,x')  = \sum_{j=0}^\infty B_j 
  (A^j b,x') ,
  \qquad x' \in \calD .
\end{equation}
\end{definition}

Thus the meaning of \eqref{eq:nud2} depends on the choice of an
admissible set $\calD$ of linear functionals.  It will be seen that in
applications natural choices are usually clear.

\begin{maintheorem}
  Let $(A,b,c)$ be given data and $\calD$ an admissible set.  Define a
  Hermitian form $\calF$ on $\calD\times \calD$ by
  \begin{equation}\label{eq:Fform}
     \calF(x',y') = \sum_{j=0}^\infty \,  \left[
    \left( A^jc,x'\right) \overline{\left( A^jc,y'\right)}
     - \left( A^jb,x'\right) \overline{\left( A^jb,y'\right)}
    \right],
   \qquad x',y'\in\calD.
  \end{equation}
  Let $\kappa$ be a nonnegative integer.
    \begin{enumerate}
    \item[(1)] If $\calF$ has $\kappa$ negative squares, there is a pair
      $(f,B)$, where $f \in \Szero$ and $B$ is a Blaschke product of
      degree~$\kappa$, such that $f(A)c = B(A)b$.
    \item[(2)] If there is a pair $(f,B)$ as in $(1)$, then $\calF$ has
      $\kappa'$ negative squares for some $\kappa' \le \kappa$.
  \end{enumerate}
\end{maintheorem}

We show that in the case $\kappa=0$ this result implies Theorem 2.3 in
\cite{RR1985}.  The reason there is something to show is that we now
require condition (iii) in Definition~\ref{D:admissible} as a
hypothesis, and this is not a hypothesis in \cite[Theorem 2.3,
p.~23]{RR1985}.  It is therefore sufficient to show that when
$\kappa=0$, condition (iii) in Definition~\ref{D:admissible} is
redundant when the other hypotheses in parts (1) and (2) of the Main
Theorem are satisfied.

(1) Assume the hypotheses of part (1) of the Main Theorem with
$\kappa=0$, but not necessarily condition (iii) in
Definition~\ref{D:admissible}.  Then $\calF$ is nonnegative, and we have
the inequality
\begin{equation}
  \notag
  \sum_{j=0}^\infty \, |\big( A^jb,x'\big)|^2
  \le
  \sum_{j=0}^\infty \, |\big( A^jc,x'\big)|^2,
  \qquad x' \in \calD,
\end{equation}
which implies (iii).

(2) Assume the hypotheses of part (2) of the Main Theorem with
$\kappa=0$, but not necessarily condition (iii) in
Definition~\ref{D:admissible}.  Then $B$ is a constant of absolute
value one, and we may take $B \equiv 1$.  Thus
\begin{equation}\label{eq:specialpsi}
  \sum_{j=0}^\infty f_j (A^j c,x')  = 
  (b,x') , \qquad x' \in \calD.
\end{equation}
If some $x_0'$ in $\calD$ annihilates $c,Ac,A^2c,\dots$, then by
Definition~\ref{D:admissible}(i), for all $k = 0,1,2,\dots$, $(A')^k
x_0'$ also belongs to $\calD$ and annihilates $c,Ac,A^2c,\dots$, and
hence $(b, (A')^kx_0') = 0$ by \eqref{eq:specialpsi}.  Thus $x_0'$
annihilates $b,Ab,A^2b,\dots$, and we again obtain
Definition~\ref{D:admissible}(iii).

The role of condition (iii) in Definition~\ref{D:admissible} can be
seen in the applications in Section~\ref{S:applications}.  For
example, in Theorem~\ref{th:pick} the condition follows from the
assumption that the points $z_1,\dots,z_n$ are distinct; this
assumption is redundant in classical Pick-Nevanlinna interpolation in
the case $\kappa=0$.

The proof of the Main Theorem is given in Section~\ref{S:proof}.

\section{Classical interpolation 
  problems on the unit disk}\label{S:applications}

We begin with the prototype for classical interpolation theory,
Pick-Nevanlinna interpolation at a finite number of points.

\begin{theorem}\label{th:pick}
  Let $z_1,\dots,z_n$ be distinct points in~$\boldD$, $w_1,\dots,w_n$
  any complex numbers, and let $\kappa$ be a nonnegative integer.  Set
  \begin{equation}
    \notag
    P = \begin{bmatrix}
    \dfrac{1 - w_j \bar w_k}{1 - z_j \bar z_k}
  \end{bmatrix}_{j,k=1}^n .
  \end{equation}
  \begin{enumerate}
  \item[(1)] If $P$ has $\kappa$ negative eigenvalues, then there is a
    pair $(f,B)$ with $f \in \Szero$ and $B$ a Blaschke product of
    degree~$\kappa$ such that $f(z_j) = B(z_j) w_j$ for all
    $j=1,\dots,n$.
  \item[(2)] If there is a pair $(f,B)$ as in $(1)$, then $P$ has
    $\kappa'\le \kappa$ negative eigenvalues.
\end{enumerate}
\end{theorem}

\begin{proof}
  We apply the Main Theorem with $\calV = \boldC^n$.  Identify
  $\calV'$ with $\boldC^n$ with the pairing $(x,y) = x_1 y_1 + \cdots
  x_n y_n$ (for any $a$ in $\boldC^n$, we write $a_j$ for its
  entries).  For the data $(A,b,c)$, choose
  \begin{equation}
  \notag
  A= \diag \{z_1,\dots,z_n\},
  \qquad 
  b = 
  \begin{bmatrix}
    w_1 \\ \vdots \\ w_n
  \end{bmatrix} ,
  \qquad 
  c = 
  \begin{bmatrix}
    1 \\ \vdots \\ 1
  \end{bmatrix} ,
  \end{equation}  
  and take $\calD = \calV'$.  The verification that $\calD$ is
  admissible is routine.  To check property (iii) in
  Definition~\ref{D:admissible}, suppose $x \in \calD$ annihilates
  $c,Ac,A^2c,\dots$.  Then 
  \begin{equation*}
    \sum_{p=1}^n z_p^j x_j = 0,\qquad  j \ge 0.
  \end{equation*}
  Since $z_1,\dots,z_n$ are distinct, $x = 0$ and trivially $x$
  annihilates $b,Ab,A^2b,\dots\;$.  A short calculation shows that the
  Hermitian form \eqref{eq:Fform} is given by
  \begin{equation}
  \notag
  \calF(x,y) = \sum_{j,k=1}^n \dfrac{1 - w_j \bar w_k}{1 - z_j \bar z_k} 
                    \;x_j \bar y_k ,
               \qquad  x,y\in\calD.
  \end{equation}
  Thus $\sqminus \calF$ is equal to the number of negative eigenvalues of
  $P$.  To see the meaning of the identity $B(A)b=f(A)c$, in
  \eqref{eq:twosums} choose the linear functional induced by the
  standard unit vector with $1$ in the $k$-th entry and all other
  entries zero.  In this case the condition reduces to $ B(z_k) w_k =
  f(z_k)$.  Thus the result is a particular case of the Main Theorem.
\end{proof}

With a slight variation of method the preceding result can be extended
to an arbitrary set of points.

\begin{theorem}\label{th:pick2}
  Let $S: \Omega \to \boldC$ be a function defined on a subset
  $\Omega$ of $\boldD$, and let $\kappa$ be a nonnegative integer.
  Set
  \begin{equation}
    \notag
    K(\zeta,z) = \frac{1-S(z)\overline{S(\zeta)}}
                 {\vphantom{X^{X^2}}1-z\bar\zeta},
    \qquad \zeta, z \in \Omega.
  \end{equation}
  \begin{enumerate}
  \item[(1)] If $\sqminus K = \kappa$, then there
    is a pair $(f,B)$ with $f \in \Szero$ and $B$ a Blaschke product
    of degree~$\kappa$ such that $f(z) = B(z) S(z)$ for all $z\in
    \Omega$. 
  \item[(2)] If there is a pair $(f,B)$ as in $(1)$, then $\sqminus K
    = \kappa' \le \kappa$.
\end{enumerate}
\end{theorem}

\begin{proof}
  Let $\calV$ be the vector space of all functions $x \colon \Omega
  \to \boldC$.  For the data $(A,b,c)$, let $A$ be multiplication by
  the independent variable on $\calV$, $b= S$, and $c \equiv 1$.  Every
  function $y$ on $\Omega$ with finite support induces a linear
  functional on $\calV$ by the formula $(x,y) = \sum_{z\in\Omega}
  x(z)y(z)$, and the set $\calD$ of all such functionals is
  admissible.  A straightforward application of the Main Theorem as in
  the proof of Theorem~\ref{th:pick} yields the result.
\end{proof}

The next result concerns the interpolation of a finite number of
derivatives, that is, it is of the Carath\'eodory-Fej\'er type.  We
use standard notation for matrices and their induced operators on
Euclidean spaces.  Write $1$ for an identity matrix or operator and
$^*$ for conjugate transpose or adjoint.

\begin{theorem}\label{th:CF}
  Let $w(z) = w_0 + w_1 z + \cdots + w_n z^n$ be a polynomial with
  complex coefficients, and set
  \begin{equation}
    \notag
    T = 
    \begin{bmatrix}
      w_0 & w_1  & \cdots &w_n  \\
      0 & w_0  &  \cdots &w_{n-1}  \\
      && \cdots && \\
      0 & 0 & \cdots & w_0
    \end{bmatrix}.
  \end{equation}
  Let $\kappa$ be a nonnegative integer such that $\kappa \le n+1$.
   \begin{enumerate}
   \item[(1)] If $1 - T^*T$ has $\kappa$ negative eigenvalues, then
     there is a pair $(f,B)$ with $f \in \Szero$ and $B$ a Blaschke
     product of degree~$\kappa$ such that $B(z)w(z) = f(z) +
     \oh(z^{n+1})$.
     
   \item[(2)] If there is a pair $(f,B)$ as in $(1)$, then $1 - T^*T$
     has $\kappa'\le \kappa$ negative eigenvalues.
  \end{enumerate}
\end{theorem}

\begin{corollary}\label{CF2}
  Let $w(z) = w_0 + w_1 z + \cdots + w_n z^n$ and $T$ be as in
  Theorem~$\ref{th:CF}$, and let~$\kappa$ be a nonnegative integer
  such that $\kappa \le n+1$.  If $1 - T^*T$ has $\kappa$ negative
  eigenvalues, there is a $\kappa'\le \kappa$ and a function $S(z)$ in
  $\Skappaprime$ which is holomorphic at the origin and such that $
  w(z) = S(z) + \oh(z^{n+1-\kappa})$.
\end{corollary}

\begin{proof}[Proof of Theorem~$\ref{th:CF}$]
  Let $\calV$ and $\calV'$ be as in the proof of Theorem~\ref{th:pick}
  but with $\boldC^n$ replaced by $\boldC^{n+1}$.  For the data
  $(A,b,c)$, choose
  \begin{equation}
  \notag
  A= 
  \begin{bmatrix}
    0&0&\cdots&0&0 \\
    1&0&\cdots&0&0 \\
    0&1&\cdots&0&0 \\
    && \cdots && \\
    0&0&\cdots&1&0
  \end{bmatrix},
  \qquad 
  b = 
  \begin{bmatrix}
    w_0 \\ w_1 \\ \vdots \\ w_n
  \end{bmatrix} ,
  \qquad 
  c = 
  \begin{bmatrix}
    1 \\ 0 \\ \vdots \\ 0
  \end{bmatrix} .
  \end{equation}  
  The set $\calD= \calV'$ is admissible.  If $\eH$ is $\boldC^{n+1}$
  in the Euclidean inner product, then $\calF(x,y) =
  {\ip{(1-T^*T)x}{y}}_{\eH}$, $x,y \in \eH$.  In fact, for any $x \in
  \eH$, 
  \begin{align*}
    \calF(x,x) &= \big( |x_0|^2 + \cdots + |x_n|^2  \big)  \\
            &\hskip1cm   - \big( |w_0 x_0 + \cdots + w_n x_n|^2 
                   + |w_0 x_1 \cdots + w_{n-1}x_n|^2 
                   + \cdots + |w_0 x_n|^2 \big) \\
    &= \left\| x \right\|_{\eH}^2 -  
                     \left\| Tx \right\|_{\eH}^2 \\[4pt]
    &= {\ip{(1-T^*T)x}{x}}_{\eH}. 
  \end{align*}
Thus $\sqminus \calF$ is equal to the number of negative eigenvalues of
$1-T^*T$.  The equation $f(A)c = B(A)b$ with $f(z) = \sum_{j=0}^\infty
f_j z^j$ and $B(z) = \sum_{j=0}^\infty B_j z^j$ is equivalent to the
identities
\begin{equation}
  \notag
   f_0 = w_0 B_0, \; 
   f_1 = w_1 B_0 + w_0 B_1,  \; 
   \dots,  \; 
   f_n = w_n B_0 + w_{n-1} B_1 + \cdots + w_0 B_n,
\end{equation}
or $B(z)w(z) = f(z) + \oh(z^{n+1})$.  The result thus follows from the
Main Theorem.
\end{proof}

\begin{proof}[Proof of Corollary~$\ref{CF2}$]
  Let $(f,B)$ be a pair as in part $(1)$ of Theorem~\ref{th:CF}.  If
  $B(z)$ has a zero of order $r$ at the origin, $f(z)$ has a zero of
  order at least $r$ at the origin.  Hence $S(z) = f(z)/B(z)$ belongs
  to $\Skappaprime$ for some $\kappa' \le \kappa$ and is holomorphic
  at the origin, and $w(z) = S(z) + \oh(z^{n-r+1})= S(z) + \oh(z^{n+1
  - \kappa})$.
\end{proof}

A simultaneous generalization of the Pick-Nevanlinna and
Carath\'eodory-Fej\'er problems can be treated in the same way by
choosing $A$ in Jordan form.  The calculations are straightforward but
somewhat lengthy, and we shall not pursue this direction.  For the
definite case, see \cite{helton} and \cite[\S2.6]{RR1985}.

The Main Theorem also yields a result on generalized interpolation in
the sense of Sarason \cite{sarason}.  Let $C$ be an inner function on
$\boldD$, and let $\eH(C) = H^2 \ominus CH^2$.  The reproducing kernel
for $\eH(C)$ is given by
\begin{equation}\label{keykernel}
K_C(w,z) = \frac{1 - C(z)\overline{C(w)}}{1 - z\bar w},
\qquad z,w \in \boldD.
\end{equation}
Let $S$ be the shift operator $S \colon h(z) \to zh(z)$ on $H^2$, and
 let $T$ be the compression of $S$ to $\eH(C)$, that is,
\begin{equation*}
T = P_{\eH(C)} S\vert_{\eH(C)},
\end{equation*}
where $P_{\eH(C)}$ is the projection operator on $H^2$ with range
$\eH(C)$.  The space $\eH(C)$ is invariant under $S^*$ and $T^* =
S^*\vert_{\eH(C)}$.  Since $T$ is completely nonunitary, for any
$\varphi \in H^\infty$ an operator $\varphi(T)$ on $\eH(C)$ is defined
by the $H^\infty$-functional calculus (see \cite{schreiber} and
\cite[p.~114]{nagyfoias}):
\begin{equation}\label{hinfinity}
\varphi(T) = s\text{-}\!\lim_{r \uparrow 1} \varphi(rT).
\end{equation}
Equivalently, for this particular situation, $\varphi(T) = P_{\eH(C)}
M_\varphi \vert_{\eH(C)}$, where $M_\varphi$ is multiplication by
$\varphi$ on $H^2$.  For every $\varphi \in H^\infty$, $\varphi(T)$
commutes with $T$, and $\varphi(T)$ is a contraction if and only if
$\varphi$ is a Schur function.

\begin{theorem}\label{generalizedinterpolation}
Let $C$ be an inner function on the unit disk, and define $T$ on
$\eH(C)$ as above.  Let $R$ be a bounded linear operator on $\eH(C)$
such that $TR=RT$.
\begin{enumerate}
\item[(1)]  If $\,1 - RR^*$ has $\kappa$ negative squares, then there is
a pair $(f,B)$, where $f \in \Szero$ and $B$ is a Blaschke product of
degree~$\kappa$, such that
\begin{equation*}
B(T)R = f(T).
\end{equation*}
\item[(2)] If there is a pair $(f,B)$ as in $(1)$, $1 - RR^*$ has
$\kappa'$ negative squares for some $\kappa' \le \kappa$.
\end{enumerate}
\end{theorem}

If $R$ is a contraction, the condition in (1) is satisfied with
$\kappa = 0$, and in this case the result reduces to the original
theorem of Sarason \cite[Theorem~1]{sarason}.

\begin{proof}
In the Main Theorem, let $\calV = \eH(C)$, $A=T$, $c= K_C(0,\cdot)$,
and $b = Rc = RK_C(0,\cdot)$.  Let $\calD$ be the set of continuous
linear functionals on $\calV = \eH(C)$; thus $\calD = \{ x_k' \colon k
\in \eH(C)\}$ where for any $k \in \eH(C)$,
\begin{equation*}
(h,x_k') = {\ip{h}{k}}_{\eH(C)},  \qquad h \in \eH(C).
\end{equation*}
Then $A'x_k' = x_{T^*k}'$ for any $k \in \eH(C)$, and so condition (i)
in Definition~\ref{D:admissible} holds.  To verify (ii), notice that
for any $k(z) = \sum_0^\infty a_j z^j$ in $\eH(C)$,
\begin{equation*}
\sum_{j=0}^\infty\,
 |(A^j c, x_k')|^2
 = \sum_{j=0}^\infty\, |{\ip{K_C(0,\cdot)}{T^{*j}k}}_{\eH(C)}|^2
 = \sum_{j=0}^\infty\, |a_j|^2 = \| k \|_{H^2}^2 < \infty.
\end{equation*}
If we replace $c$ by $b$ and use the identity $RT=TR$, we obtain
\begin{equation*}
\sum_{j=0}^\infty\, |(A^j b, x_k')|^2
 = \sum_{j=0}^\infty\, |{\ip{K_C(0,\cdot)}{T^{*j}R^*k}}_{\eH(C)}|^2
 = \|R^* k \|_{H^2}^2 < \infty,
\end{equation*}
and thus condition (ii) in Definition~\ref{D:admissible} holds.
Condition (iii) holds trivially because the only $k \in \eH(C)$ such
that 
\begin{equation*}
  (A^j c, x_k') = {\ip{K_C(0,\cdot)}{T^{*j}k}}_{\eH(C)} = 0
\end{equation*}
for all $j \ge 0$ is $k=0$.  We have shown that $\calD$ is admissible.
The form \eqref{eq:Fform} is given by
\begin{align}
\calF(x_h',x_k') &=
   \sum_{j=0}^\infty\, \Big[ 
 {\ip{T^jK_C(0,\cdot)}{h}}_{\eH(C)} {\ip{k}{T^jK_C(0,\cdot)}}_{\eH(C)}
     \label{RRstar} \\
&\hskip2cm 
- {\ip{T^jRK_C(0,\cdot)}{h}}_{\eH(C)} {\ip{k}{T^jRK_C(0,\cdot)}}_{\eH(C)}
 \Big] \notag \\[4pt]
   &= {\ip{k}{h}}_{H^2} - {\ip{R^*k}{R^*h}}_{H^2} \notag \\[4pt]
   &=  {\ip{k}{h}}_{\eH(C)} - {\ip{R^*k}{R^*h}}_{\eH(C)} \notag \\[4pt]
   &= {\ip{(1 - RR^*)k}{h}}_{\eH(C)}
\notag
\end{align}
for any $h$ and $k$ in $\eH(C)$.  

(1) Assume that $1-RR^*$ has $\kappa$ negative squares.  By
\eqref{RRstar}, the Hermitian form \eqref{eq:Fform} has $\kappa$
negative squares.  Hence by part (1) of the Main Theorem, there is a
function $f \in \Szero$ and a Blaschke product $B$ of degree~$\kappa$
such that $B(A)b = f(A)c$ in the sense of Definition~\ref{nudelmaneq},
that is, if $B(z) = \sum_0^\infty B_j z^j$ and $f(z) = \sum_0^\infty
f_j z^j$, then for every $h \in \eH(C)$,
\begin{equation*}
\sum_{j=0}^\infty \,
B_j{\ip{T^jRK_C(0,\cdot)}{h}}_{\eH(C)}
=
\sum_{j=0}^\infty \,
f_j{\ip{T^j K_C(0,\cdot)}{h}}_{\eH(C)}.
\end{equation*}
Using Abel summation of these series, we see that 
\begin{equation}\label{abel}
B(T)RK_C(0,\cdot) = f(T)K_C(0,\cdot).
\end{equation}
Since $R$ commutes with $T$, it commutes with $B(T)$ and $f(T)$.
Hence $B(T)R$ and $f(T)$ agree on the smallest invariant subspace of
$T$ containing $K_C(0,\cdot)$.  The latter subspace is all of
$\eH(C)$, and thus we obtain $B(T)R = f(T)$.

(2) Assume that a pair $(f,B)$ exists as in (1).  Reversing the
    preceding steps, we see that $B(A)b=f(A)c$, hence by part (2) of
    the Main Theorem the form \eqref{eq:Fform} has $\kappa$ negative
    squares.  Therefore by \eqref{RRstar}, $1 - RR^*$ has $\kappa'$
    negative squares for some $\kappa' \le \kappa$.
\end{proof}

More generally, let $C$ be a function which is holomorphic and bounded
by one on $\boldD$ (not necessarily an inner function), and let
$\eH(C)$ be the Hilbert space with reproducing kernel
\eqref{keykernel} \cite{dbr, sarasonbook}.  Then $\eH(C)$ is contained
contractively in $H^2$, that is, the inclusion mapping $E \colon
\eH(C) \to H^2$ is a contraction operator (the inclusion is isometric
if and only if $C$ is an inner function).  There is a contraction
operator $T$ on $\eH(C)$ such that $T^*:h(z) \to [h(z)-h(0)]/z$ for
every $h(z)$ in $\eH(C)$.  In view of the difference-quotient
inequality
\begin{equation*}
\| [h(z)-h(0)]/z \|_{\eH(C)}^2 \le \| h(z) \|_{\eH(C)}^2 - |h(0)|^2,
\end{equation*}
which holds for every element $h(z)$ of $\eH(C)$, $T$ is completely
nonunitary.  Hence for any $\varphi \in H^\infty$ we may define
$\varphi(T)$ by \eqref{hinfinity}.  The selfadjoint operators $G =
E^*E$ on $\eH(C)$ and $D = EE^*$ on $H^2$ also play a role.  They
satisfy $0 \le G \le 1$ and $0 \le D \le 1$.  It is not hard to see
that
\begin{equation*}
D = 1 - M_C M_C^*,
\end{equation*}
where $M_C$ is multiplication by $C$ on $H^2$.  Let
\begin{equation*}
 \eM = H^2 \ominus \ker\,D.
\end{equation*}
Then $\eH(C)$ coincides with the range of $D^{1/2}$ in the ``range
norm'', that is, the unique norm such that $D^{1/2}$ acts as a partial
isometry from $H^2$ onto $\eH(C)$ with initial space $\eM$.  One can
use the relation $DE=EG$ to show that $G$ is unitarily equivalent to
$D \vert_\eM$ by means of the natural isomorphism $U =
D^{1/2}\vert_\eM$ from $\eM$ onto $\eH(C)$.

\begin{theorem}\label{gengeninterp}
  Let $C$ be a function which is holomorphic and bounded by one on the
  unit disk, and let $T$ be the operator on $\eH(C)$ such that
  $T^*:h(z) \to [h(z)-h(0)]/z$ for every $h(z)$ in $\eH(C)$.  Let
  $G=E^*E$ and $E$ be as above.  Let $R$ be a bounded linear operator
  on $\eH(C)$ such that $TR=RT$.
\begin{enumerate}
\item[(1)]  If $G - RGR^*$ has $\kappa$ negative squares, then there is
a pair $(f,B)$, where $f \in \Szero$ and $B$ is a Blaschke product of
degree~$\kappa$, such that
\begin{equation*}
B(T)R = f(T).
\end{equation*}
\item[(2)] If there is a pair $(f,B)$ as in $(1)$, $G - RGR^*$ has
$\kappa'$ negative squares for some $\kappa' \le \kappa$.
\end{enumerate}
\end{theorem}

\noindent {\it Remarks}. (1)  Let $f \in \Szero$, and define $R=f(T)$ by
the $H^\infty$-functional calculus.  Set $\tilde f(z) =
\overline{f(\bar z)}$.  Then $RT=TR$.  Since
$Ef(T)^* = E\tilde f(T^*) = \tilde f(S^*)E$, we obtain
\begin{equation*}
G - RGR^* = E^*E - f(T)E^*Ef(T)^* = E^*E - E^*\tilde f(S^*)^*
\tilde f(S^*)E \ge 0,
\end{equation*}
because $\tilde f(S^*)$ is a contraction.

(2)  It is a corollary of the theorem that if $R$ is a
bounded operator on $\eH(C)$ such that $RT=TR$ and $G-RGR^* \ge 0$,
then $1 - RR^* \ge 0$.  Indeed, by the theorem $R=f(T)$ for some $f
\in \Szero$, and hence $R$ is a contraction.

\begin{proof}[Proof of Theorem~$\ref{gengeninterp}$]
We proceed as in the proof of Theorem~\ref{generalizedinterpolation}
with the same choice of $\calV$, data $(A,b,c)$, and set $\calD$ of
linear functionals.  In exactly the same way, we
show that $\calD$ is admissible.  The calculation of \eqref{eq:Fform}
is the same except for the last step.  Now we have
\begin{equation*}
\calF(x_h',x_k')
   = {\ip{Ek}{Eh}}_{H^2} - {\ip{ERk}{ERh}}_{H^2}
   = {\ip{(G - RGR^*)k}{h}}_{\eH(C)}
\end{equation*}
for all $h$ and $k$ in $\eH(C)$.

(1) If $G-RGR^*$ has $\kappa$ negative squares, so does
\eqref{eq:Fform}., and hence by the Main Theorem there is a pair
$(f,B)$, where $f \in \Szero$ and $B$ is a Blaschke product of
degree~$\kappa$, such that \eqref{eq:nud2} holds.  We use Abel
summation as in the proof of Theorem~\ref{generalizedinterpolation} to
show that this implies \eqref{abel}, where $B(T)$ and $f(T)$ are
defined by the $H^\infty$-functional calculus.  The commutivity of $T$
and $R$ then implies that $B(T)R = f(T)$.

(2) Reverse the steps in (1).
\end{proof}

We turn now to problems of Loewner type, which concern boundary data.
In the scalar case, solutions of interpolation problems of the type
that we treat here are unique, and the results that we obtain are
characterizations of boundary functions and their restrictions.

Let $\sigma$ be normalized Lebesgue measure on $\partial\boldD$:
$|u|=1$.  When no confusion can arise, we write $L^2,L^\infty$ for
$L^2(\sigma),L^\infty(\sigma)$; we identify $H^2$ with a subspace of
$L^2$ in the usual way.  If $\Delta$ is a measurable subset of
$\partial\boldD$, $L^2(\Delta)$ is the subspace of functions in $L^2$
supported on $\Delta$.  The characteristic function of $\Delta$ is
denoted $1_\Delta$; if $\varphi$ is a function on $\Delta$, we view
$\varphi1_\Delta$ as defined on all of $\partial\boldD$ and equal to
zero on the complement of $\Delta$.

We note some preliminary results from function theory.

\begin{enumerate}
\item[(i)] If $\varphi, \psi \in L^2(\Delta)$, then
  \cite[p.~31]{RR1985}
\begin{multline}\label{keylimit}
   \lim_{r\uparrow 1} \int_\Delta \int_\Delta 
   \frac{\varphi(u)\overline{\psi(v)}}{1 - r^2 u \bar v}
        \; d\sigma(u)d\sigma(v)   \\
 =  \sum_{j=0}^\infty 
     \Big( \int_\Delta u^j \varphi(u) \; d\sigma(u)  \Big)
        \Big( \int_\Delta v^j \psi(v) \; d\sigma(v)  \Big)^{-\!\!-}     
   = {\ip{Q_-(\varphi1_\Delta)}{\psi1_\Delta}}_{L^2},
\end{multline}
where $Q_-$ is the orthogonal projection on $L^2$ whose range is the
closed span $L^2_-$ of all functions $u^j$, $j \le 0$.  In particular,
the limit on the left exists.
\item[(ii)] If $B$ is an inner function and $M_B$ is multiplication by
  $B$ on $L^2$, then
\begin{equation}\label{hbprojection} 
M_BQ_-M_B^* = Q_- + P_{u\eH(B)},
\end{equation} 
where $\eH(B) = H^2 \ominus BH^2$ and $P_{u\eH(B)}$ is the orthogonal
projection on $L^2$ whose range is $u\eH(B)$.  To prove this, write
$L^2 = L^2_- \oplus uBH^2 \oplus u\eH(B)$.  Both sides of
\eqref{hbprojection} coincide with the identity operator on $L^2_-$,
and both coincide with the zero operator on $uBH^2$.  It remains to
check the actions of each side on elements of $L^2$ of the form $uh$
with $h$ in $\eH(B)$.  This amounts to showing that $BQ_-(\bar B uh) =
uh$, or equivalently, that $Q_-(\bar B uh) = \bar B uh$; the last
equation holds by the characterization of $\eH(B)$ as the set of
elements $h$ of $H^2$ such that $\bar B h \perp H^2$.

\item[(iii)] (Kronecker's Theorem) This result characterizes
  finite-rank Hankel operators and has several formulations,
  including: {\it Let $F \in L^\infty$, and assume that the operator
    $H_F$ on $H^2$ to $L^2 \ominus H^2$ which is defined as
    compression of multiplication by $F$ has rank~$\kappa$.  Then
    there is a Blaschke product $B$ of degree $\kappa$ such that $BF
    \in H^\infty$.}  For example, see Peller \cite[p.~77]{peller}.
  Here $H_F h = P_- (Fh)$, $h \in H^2$, where $P_-$ is the projection
  onto $L^2 \ominus H^2$, equivalently, $H_F h = \bar u Q_- (uFh)$
  where $Q_-$ is as above.
\end{enumerate}

Our first boundary result characterizes boundary functions of
functions in $\Skappa$. 

\begin{theorem}\label{th:boundaryD}
  Let $b$ and $c$ be complex-valued measurable functions on a Borel
  set $\Delta \subseteq \partial\boldD$ such that $c \neq 0$
  $\sigma$-a.e.\ on $\Delta$.  Let $\calD$ be the set of functions
  $\varphi$ on $\Delta$ such that $b\varphi,c\varphi\in L^2(\Delta)$.
  If $\kappa$ is a nonnegative integer, then there is a function
  $S(z)$ in $\Skappa$ whose boundary function satisfies $b(u) =
  S(u)c(u)$ $\sigma$-a.e.\ on $\Delta$ if and only if the Hermitian
  form
\begin{equation}\label{eq:bcform}
L(\varphi,\psi) = \lim_{r\uparrow 1} \int_\Delta \int_\Delta 
 \frac{c(u)\overline{c(v)} - b(u)\overline{b(v)}}{1 - r^2 u \bar v}
  \,\varphi(u)\overline{\psi(v)}\; d\sigma(u)d\sigma(v),
\end{equation} 
$\varphi,\psi\in\calD$, has $\kappa$ negative squares.
\end{theorem}

\begin{proof}
  Case 1: $\partial\boldD \setminus \Delta$ is not a $\sigma$-null
  set.

  We apply the Main Theorem with $\calV = \{ pb+qc : p,q\;
  \text{polynomials}\}$.  For the data $(A,b,c)$, let $A$ be
  multiplication by the independent variable, and let $b$ and $c$ be
  the given functions.  We view $\calD$ as a subspace of $\calV'$ by
  identifying each $\varphi \in \calD$ with a linear functional on
  $\calV$ by writing
  $$
  (f,\varphi) = \int_\Delta f\varphi \; d\sigma,\qquad
  f\in\calV.
  $$
  Condition (i) in the definition of an admissible family is easily
  verified, and (ii) holds because the Fourier coefficients of a
  function in $L^2$ are square summable.  To check (iii),
  suppose $\varphi \in \calD$ annihilates $c,Ac,A^2c,\dots$.  Then
  \begin{equation}
    \notag
    \int_{\partial\boldD} u^j c \varphi 1_\Delta \; d\sigma
    =0,
    \qquad j = 0,1,2,\dots.
  \end{equation}
  Thus $c \varphi 1_\Delta$ belongs to $u H^2$.  Since
  $\sigma(\partial\boldD \setminus \Delta) > 0$, $c \varphi 1_\Delta =
  0$ $\sigma$-a.e.\ on $\partial\boldD$ and hence $\varphi = 0$
  $\sigma$-a.e.\ on $\Delta$.  Therefore $\varphi$ is the zero element
  of $\calD$ and hence annihilates $b,Ab,A^2b,\dots$.  Thus $\calD$ is
  admissible.
  
  By \eqref{keylimit}, the Hermitian form \eqref{eq:Fform} coincides
  with \eqref{eq:bcform}.  An argument in \cite[p.~31]{RR1985} shows
  that for any Schur functions $f$ and $B$, the identity $B(A)b =
  f(A)c$ is equivalent to the relation $Bb=fc$
  $\sigma$-a.e.\ on~$\Delta$.
  
  Now suppose $\sqminus L = \kappa$.  By part (1) of the Main Theorem,
  there is a pair $(f,B)$, $f \in \Szero$ and $B$ a Blaschke product
  of degree~$\kappa$, such that $f(A)c=B(A)b$ and hence $Bb=fc$
  $\sigma$-a.e.\ on~$\Delta$.  Hence $b = Sc$
  $\sigma$-a.e.\ on~$\Delta$, where $S = f/B$ is in $\Skappaprime$ for
  some $\kappa' \le \kappa$.  Let $S = f_1/B_1$ be a Kre{\u\i}n-Langer
  factorization of $S$, so $B_1$ is a Blaschke product of
  degree~$\kappa'$ (which is obtained by cancelling all of the common
  simple Blaschke factors in $f$ and $B$).  Then $f_1c = B_1b$
  $\sigma$-a.e.\ on~$\Delta$ and so $f_1(A)c = B_1(A)b$.  Part (2) of
  the Main Theorem then gives $\sqminus L \le \kappa'$.  Therefore
  $\kappa = \kappa'$ and the sufficiency part of the theorem follows.
 
  Conversely, let $b=Sc$ $\sigma$-a.e.\ on~$\Delta$ where $S \in
  \Skappa$.  If $S = f/B$ is a Kre{\u\i}n-Langer factorization of $S$,
  then $f(A)c=B(A)b$, and by part (2) of the Main Theorem, $\sqminus L
  = \kappa' \le \kappa$.  By the proof of sufficiency above, $b =
  S_1c$ $\sigma$-a.e.\ on $\Delta$, where $S_1 \in \Skappaprime$.
  Since we assume that $c \neq 0$ $\sigma$-a.e.\ on $\Delta$, $S_1 =
  b/c = S \in \Skappa$, and hence $L$ has $\kappa' = \kappa$ negative
  squares.  The necessity part of the theorem follows.

\noindent  Case 2: $\Delta = \partial\boldD$.

  Given $\varepsilon > 0$, let $\Delta_\varepsilon$ be the set of
  points $u \in \partial\boldD$ with $\varepsilon \le \arg \, u \le
  2\pi$.  Let $b_\varepsilon,c_\varepsilon$ be the restrictions of
  $b,c$ to $\Delta_\varepsilon$.  Let $\calD_\varepsilon$ be the set
  of restrictions of functions in $\calD$ to $\Delta_\varepsilon$, and
  define $L_\varepsilon$ as the restriction of $L$ to
  $\calD_\varepsilon \times \calD_\varepsilon$.

  If $L$ has $\kappa$ negative squares, so does $L_\varepsilon$ for
  all sufficiently small $\varepsilon$.  By Case~1, for such
  $\varepsilon$ there is a function $S_\varepsilon(z)$ in $\Skappa$
  whose boundary function satisfies $b(u) = S_\varepsilon(u)c(u)$
  $\sigma$-a.e.\ on $\Delta_\varepsilon$.  Since functions in $\Skappa$
  are determined by their boundary values on a set of positive
  measure, $S_\varepsilon = S$ is independent of~$\varepsilon$.  By
  construction, $b(u) = S(u)c(u)$ $\sigma$-a.e.\ on~$\partial\boldD$.

  In the other direction, if there is an $S$ in $\Skappa$ whose
  boundary function satisfies $b(u) = S(u)c(u)$
  $\sigma$-a.e.\ on~$\partial\boldD$, the same relation holds
  $\sigma$-a.e.\ on~$\Delta_\varepsilon$.  Again by Case~1,
  $L_\varepsilon$ has $\kappa$ negative squares for all sufficiently
  small $\varepsilon$.  Hence $L$ has $\kappa$ negative squares.  Thus
  the result holds in the case $\Delta=\partial\boldD$.
\end{proof}

\begin{corollary}\label{kronecker}
   Let $S_0$ be a complex-valued measurable function on
  $\partial\boldD$, and let $\kappa$ be a nonnegative integer.  Then
  there is a function $S(z)$ in $\Skappa$ whose boundary function
  satisfies $S(u) = S_0(u)$ $\sigma$-a.e.\ on $\partial\boldD$ if and
  only if the Hermitian form
\begin{equation}\label{kroneckerform}
L(\varphi,\psi) = \lim_{r\uparrow 1} \int_{\partial\boldD}
  \int_{\partial\boldD}
 \frac{1 - S_0(u)\overline{S_0(v)}}{1 - r^2 u \bar v}
  \,\varphi(u)\overline{\psi(v)}\; d\sigma(u)d\sigma(v),
\end{equation} 
$\varphi, S_0\varphi,\psi,S_0\psi\in L^2$, has $\kappa$ negative
squares.  When this condition is satisfied, then $|S_0(u)| \le 1$
$\sigma$-a.e.\ on $\partial\boldD$.
\end{corollary}

\begin{proof}
In Theorem~\ref{th:boundaryD} choose $\Delta = \partial\boldD$, $c=1$,
and $b=S_0$.  The last statement follows by the Kre{\u\i}n-Langer
representation \eqref{eq:kreinlanger} of a generalized Schur function.   
\end{proof}

We give a different version of the corollary that assumes $|S_0(u)|
\le 1$ $\sigma$-a.e.\ on $\partial\boldD$ as a hypothesis and then
obtains additional information.  Namely, the result remains true if we
replace the form $L(\varphi,\psi)$ defined on $L^2 \times L^2$ by
\eqref{kroneckerform} by its restriction \eqref{kroneckerform2} to
$uH^2 \times uH^2$.  \medskip

\noindent {\bf Alternative form of Corollary~${\bf \ref{kronecker}}$}
{\it Let $S_0$ be a complex-valued measurable function on
   $\partial\boldD$ satisfying $|S_0(u)| \le 1$ $\sigma$-a.e.\ on
   $\partial\boldD$, and let $\kappa$ be a nonnegative integer.  Then
   there is a function $S(z)$ in $\Skappa$ whose boundary function
   satisfies $S(u) = S_0(u)$ $\sigma$-a.e.\ on $\partial\boldD$ if and
   only if the Hermitian form
\begin{equation}\label{kroneckerform2}
L_+(\varphi,\psi) = \lim_{r\uparrow 1} \int_{\partial\boldD}
  \int_{\partial\boldD}
 \frac{1 - S_0(u)\overline{S_0(v)}}{1 - r^2 u \bar v}
  \,\varphi(u)\overline{\psi(v)}\; d\sigma(u)d\sigma(v),
\end{equation} 
$\varphi,\psi\in uH^2$, has $\kappa$ negative squares.
}

\begin{proof}
  First assume that $S_0 \in \Skappa$, or, more precisely, that $S_0(u)
  = S(u)$ $\sigma$-a.e.\ for some $S(z)$ in $\Skappa$.  We show that
  $\sqminus L_+ = \kappa'$ for some $\kappa' \le \kappa$.  Let $S =
  f/B$ be a Kre{\u\i}n-Langer factorization.  Then by \eqref{keylimit}
  and \eqref{hbprojection}, for any $\varphi,\psi \in uH^2$,
\begin{align} 
L_+(\varphi,\psi) &=
  {\ip{Q_-\varphi}{\psi}}_{L^2} - {\ip{Q_- (\bar Bf\varphi)}{\bar Bf
  \psi}}_{L^2} \notag \\[4pt] 
&= -{\ip{M_B Q_- M_B^* f\varphi}{f \psi}}_{L^2} 
  \notag \\
&= -{\ip{\big(Q_- +P_{u\eH(B)}\big) f\varphi}{f \psi}}_{L^2}  \notag \\
&= -{\ip{P_{u\eH(B)} f\varphi}{f \psi}}_{L^2} . \notag
\end{align} 
Since $\eH(B)$ has dimension $\kappa$, $\sqminus L_+ = \kappa'$ for
some $\kappa' \le \kappa$.

Next suppose that $\sqminus L_+ = \kappa'$ for some $\kappa'$.  As
above, for any $\varphi = uh$ and $\psi = uk$ in $uH^2$, 
\begin{equation*} 
L_+(\varphi,\psi) =
  {\ip{Q_-\varphi}{\psi}}_{L^2} - {\ip{Q_- S_0\varphi}{S_0
  \psi}}_{L^2}
= -{\ip{Q_- uS_0 h}{uS_0 k}}_{L^2}  .
\end{equation*} 
Since we assume that $\sqminus L_+ = \kappa'$, the rank of the Hankel
operator defined by $H_{S_0}h = \bar u Q_- (uS_0 h)$, $h \in H^2$, is
$\kappa'$.  By Kronecker's theorem as stated above, there is a
Blaschke product $B$ of degree $\kappa'$ such that $BS_0 = f$ where $f
\in \Szero$ (recall that we assume that $|S_0(u)| \le 1$
$\sigma$-a.e.\ on $\partial\boldD$).  It follows that $S_0 \in
\Skappadoubleprime$ for some $\kappa'' \le \kappa'$, that is, $S_0$ is
the boundary function of a function in $\Skappadoubleprime$.  But by
the first part of the proof, we then have $\kappa' = \sqminus L_+ \le
\kappa''$, and so $\kappa'' = \kappa'$ and $S_0 \in \Skappaprime$.

The result follows on combining the two parts of the argument.
\end{proof}

\noindent {\bf Example.}  In most cases, condition (iii) in the
definition of an admissible family is trivially satisfied, but there
are situations in which it does not hold.  An example arises in the
proof of Theorem~\ref{th:boundaryD} when $\Delta = \partial\boldD$; it
is not possible to include this case in the main argument there
because condition (iii) for an admissible family may fail.  To see
this, in the proof of Theorem~\ref{th:boundaryD} (Case~1) allow
$\Delta = \partial\boldD$, and take $c = 1$ and $b = \bar u^\kappa$
for some positive integer $\kappa$.  An element $\varphi$ of $\calD =
L^2$ annihilates $c,Ac,A^2c,\dots$ if and only if
\begin{equation}
  \notag
  (A^jc,\varphi) = \int_{\partial\boldD} u^j \varphi \; d\sigma = 0,
  \qquad j \ge 0.
\end{equation}
This implies that
\begin{equation}
  \notag
  (A^jb,\varphi) = \int_{\partial\boldD} u^{j-\kappa} \varphi \; 
d\sigma = 0,
  \qquad j \ge \kappa,
\end{equation}
but the last identity can fail for $0 \le j < \kappa$, that is,
condition (iii) in the definition of an admissible family does not
hold.  Nevertheless, the Hermitian form \eqref{eq:bcform} has $\kappa$
negative squares; indeed by \eqref{keylimit} for $\varphi,\psi$ in
$\calD$,
\begin{align*}
L(\varphi,\psi) &= \lim_{r\uparrow 1} \int_\Delta \int_\Delta
 \frac{c(u)\overline{c(v)} - b(u)\overline{b(v)}}{1 - r^2 u \bar v}
 \,\varphi(u)\overline{\psi(v)}\; d\sigma(u)d\sigma(v)  \\[4pt]
&= - {\ip{\big( M_{u^\kappa} Q_- M_{u^\kappa}^* - 
Q_-\big)\varphi}{\psi}}_{L^2},
\end{align*} 
where $M_{u^\kappa} Q_- M_{u^\kappa}^* - Q_-$ is the projection of
$L^2$ onto the span of $1,u,\dots,u^\kappa$.  Of course, even in this
case the conclusion of Theorem~\ref{th:boundaryD} is true as explained
in Case~2 of the proof.  Indeed, the identity $b(u) = S(u)c(u)$ holds
$\sigma$-a.e.\ on $\partial\boldD$ where $S(z) = 1/z^\kappa$ belongs to
$\Skappa$.  \medskip

\section{Interior and boundary interpolation on a
 half-plane}\label{S:halfplane}

Half-plane results are more simply derived by direct application of
the Main Theorem than by change of variables from the disk case.  A
Blaschke product of degree $\kappa$ on the upper half-plane $\boldC_+$
is a function of the form
\begin{equation}
  \notag
 B(z) = c \prod_{j=1}^\kappa\; \frac{z-a_j}{z - \bar a_j},
\end{equation}
where $c$ is a constant of modulus one and $a_1,\dots,a_\kappa$ are
points in $\boldC_+$ (not necessarily distinct).  Every such function
can be written $B(z) = B_0((z-i)/(z+i))$, where $B_0(z)$ is a Blaschke
product of degree $\kappa$ on the unit disk, and conversely
\cite{RR1994}.  We say that a meromorphic function $S(z)$ on
$\boldC_+$ belongs to $\SkappaCplus$ if it has the form
$$S(z) = S_0 \Big( \frac{z-i}{z+i} \Big),
$$
where $S_0(z)$ belongs to $\Skappa$ as a function on $\boldD$.  If
$S(z)$ belongs to $\SkappaCplus$, then
\begin{equation}\label{linearfrac}
  f(z) = i\, \frac{1 + S(z)}{1 - S(z)}
\end{equation}
defines a function in $\Nkappa$, and every function in $\Nkappa$ is
obtained in this way; when $\kappa=0$ we exclude $S(z)\equiv 1$ from
this correspondence.

If $\Delta$ is a (Lebesgue) measurable subset of $(-\infty,\infty)$
and $\varphi, \psi \in L^2(\Delta)$, then \cite[pp.~33--34]{RR1985}
\begin{align}
   \lim_{\epsilon \downarrow 0 }\,
  \frac{i}{2}& \int_\Delta \int_\Delta
   \frac{\varphi(s)\overline{\psi(t)}}{s-t+i\epsilon}\; ds \, dt
  \label{keylimit2} \\
    &=
  \sum_{j=0}^\infty \left(
    \int_\Delta \bigg( \frac{t-i}{t+i}  \bigg)^j \frac{1}{t+i}
   \, \varphi(t) \; dt
  \right)
  \left(
    \int_\Delta \bigg( \frac{t-i}{t+i}  \bigg)^j \frac{1}{t+i}
   \, \psi(t) \; dt
  \right)^{\!-\!\!-} \notag \\[4pt]
  &= \pi\, {\ip{Q_-(\varphi1_\Delta)}{\psi1_\Delta}}_{L^2(-\infty,\infty)}
  \, . 
  \notag
\end{align}
The meaning of $Q_-$ here is different from that of \eqref{keylimit}.
Let $H^2$ now denote the Hardy class for the upper half-plane
\cite{RR1994}.  We may alternatively view $H^2$ as a subspace of
$L^2(-\infty,\infty)$ by passing to boundary values.  In
\eqref{keylimit2}, $Q_-$ is the projection of $L^2(-\infty,\infty)$
onto the orthogonal complement of $H^2$ in $L^2(-\infty,\infty)$, or
equivalently onto the set of complex conjugates of functions in $H^2$,
since
\begin{equation}
  \notag
  L^2(-\infty,\infty) = H^2 \oplus \overline{H^2} .
\end{equation}
The last equality in
\eqref{keylimit2} follows from the fact that the functions
\begin{equation}\label{basis}
  \frac{1}{\sqrt{\pi}}\bigg( \frac{t-i}{t+i}  \bigg)^j \frac{1}{t+i},
  \qquad j = 0, \pm 1, \pm 2, \dots ,
\end{equation}
are an orthonormal basis for $L^2(-\infty,\infty)$ and the part of
\eqref{basis} for $j \ge 0$ is an orthonormal basis for $H^2$.  In
particular, the infinite series in the middle term of
\eqref{keylimit2} always converges, and the limit on the left side of
\eqref{keylimit2} always exists.

The next result is a half-plane counterpart of
Theorem~\ref{th:boundaryD} and preliminary to boundary theorems for
the generalized Nevanlinna class.

\begin{theorem}\label{th:boundaryR}
  Let $b$ and $c$ be complex-valued measurable functions on a Borel
  subset $\Delta$ of $(-\infty,\infty)$ such that $c \neq 0$ a.e.\ on
  $\Delta$.  Let $\calD$ be the set of measurable functions $\varphi$
  on $\Delta$ such that $b\varphi,c\varphi\in L^2(\Delta)$.  If
  $\kappa$ is a nonnegative integer, then there is a function $S(z)$
  in $\SkappaCplus$ whose boundary function satisfies $b(x)=S(x)c(x)$
  a.e.\ on $\Delta$ if and only if the Hermitian form
\begin{equation}\label{eq:bcform2}
  L(\varphi,\psi) = \lim_{\varepsilon \downarrow 0} \; \frac{i}{2} 
      \int_\Delta \int_\Delta 
 \frac{ c(s)\overline{c(t)} - b(s)\overline{b(t)}}{s-t+i\varepsilon}
  \,\varphi(s)\overline{\psi(t)}\; ds\,dt, 
\qquad \varphi,\psi \in \calD,
\end{equation}
has $\kappa$ negative squares.
\end{theorem}

\begin{proof}
  Case 1: $(-\infty,\infty) \setminus \Delta$ is not a Lebesgue null
  set.

  Let $\calV$ be the vector space of complex-valued functions on
  $\Delta$ of the form 
  \begin{equation}
    \notag
    p\left(\dfrac{x-i}{x+i}\right) b(x) 
      + q\left(\dfrac{x-i}{x+i}\right) c(x) ,
    \qquad p,q\; \text{polynomials.}
  \end{equation}
  We identify any $\varphi \in \calD$ with the linear functional on
  $\calV$ defined by
\begin{equation}
  \notag
  (h, \varphi) = \int_\Delta \dfrac{h(t) \varphi(t)}{t+i} \; dt ,
  \qquad h \in \calV.
\end{equation}
Define data $(A,b,c)$ by choosing $A$ to be multiplication by
$(x-i)/(x+i)$ and letting $b$ and $c$ be the given functions.
Arguments as in the proof of Theorem~\ref{th:boundaryD} show that
$\calD$ is admissible; the assumption that $(-\infty,\infty) \setminus
\Delta$ is not a Lebesgue null set is used to verify condition (iii)
for an admissible set.

By \eqref{keylimit2}, the Hermitian form \eqref{eq:Fform} coincides
with \eqref{eq:bcform2}.  We show that $f_0(A)c = B_0(A)b$ for a pair
$(f_0,B_0)$ of holomorphic functions which are bounded by one on the
unit disk $\boldD$ if and only if $B(x)b(x) = f(x)c(x)$ a.e.\ on
$\Delta$, where $f$ and $B$ are defined on the upper half-plane
$\boldC_+$ by
\begin{align}
  f(z) &= f_0 \left( \dfrac{z-i}{z+i} \right) = \sum_0^\infty\, f_j
  \left( \dfrac{z-i}{z+i}  \right)^j , \notag \\
  B(z) &= B_0 \left( \dfrac{z-i}{z+i} \right) = \sum_0^\infty \,
  B_j \left( \dfrac{z-i}{z+i} \right)^j .\notag
\end{align}
If $f_0(A)c = B_0(A)b$, then for every 
 $\varphi
\in \calD$,  
$$\sum_0^\infty \, f_j (A^jc,\varphi) = \sum_0^\infty \, B_j
(A^jb,\varphi),$$
that is,
\begin{equation}
  \notag
  \sum_0^\infty \, f_j 
     \int_\Delta \dfrac{1}{t+i}\left( \dfrac{t-i}{t+i} \right)^j c(t)
  \varphi(t)\; dt
     = 
  \sum_0^\infty \, B_j \int_\Delta \dfrac{1}{t+i}\left( \dfrac{t-i}{t+i} \right)^j b(t)
  \varphi(t)\; dt.
\end{equation}
By the arbitrariness of $\varphi$, on summing the series in the Abel
sense as in \cite[p.~32]{RR1985}, we obtain 
\begin{equation}
  \notag
  f_0 \left( \dfrac{x-i}{x+i} \right)c(x) = B_0 \left(
  \dfrac{x-i}{x+i} \right)b(x),
\end{equation}
and hence $f(x)c(x)=B(x) b(x)$ a.e.\ on $\Delta$.  These steps are
reversible, and the assertion follows.

To complete the proof in Case~1, we apply the Main Theorem as in the
proof of Theorem~\ref{th:boundaryD}.

\noindent  Case 2: $\Delta = (-\infty,\infty)$.

We reduce this to Case~1 in the same way as in the proof of
Theorem~\ref{th:boundaryD}.
\end{proof}

An indefinite extension of the Generalized Loewner Theorem
\cite[p.~34]{RR1985} follows as an immediate consequence of
Theorem~\ref{th:boundaryR}.

\begin{theorem}\label{th:boundaryL}
  Let $f_0$ be a complex-valued measurable function defined on a Borel
  subset $\Delta$ of $(-\infty,\infty)$ such that $f_0 + i \neq 0$
  a.e.\ on $\Delta$.  Let $\calD$ be the set of measurable functions
  $\varphi$ on $\Delta$ such that $\varphi, f_0 \varphi \in
  L^2(\Delta)$.  If $\kappa$ is a nonnegative integer, there exists a
  function $f\in\Nkappa$ such that $f=f_0$ a.e.\ on $\Delta$ if and
  only if the Hermitian form
\begin{equation}\label{eq:bcform3}
L(\varphi,\psi) = \lim_{\varepsilon\downarrow 0} \int_\Delta \int_\Delta 
 \frac{f_0(s) - \overline{f_0(t)}}{s-t+i\varepsilon}
  \,\varphi(s)\overline{\psi(t)}\; ds\,dt,
\qquad \varphi,\psi\in\calD,
\end{equation}
has $\kappa$ negative squares.
\end{theorem}

The assumption that $f_0 + i \neq 0$ a.e.\ on $\Delta$ does not
restrict the generality: it is automatically satisfied if $f=f_0$
a.e.\ on $\Delta$ for some $f\in\Nkappa$.  For if $f(z)$ and $S(z)$ are
related by \eqref{linearfrac}, then $|S(x)| \le 1$ by the
Kre{\u\i}n-Langer representation \eqref{eq:kreinlanger}, and hence
$\text{Im}\, f_0(x) \ge 0$ a.e.\ on $\Delta$.

\begin{proof}
  We apply Theorem~\ref{th:boundaryR} with $b= f_0-i$ and $c = f_0 +
  i$.  With this choice of $b$ and $c$, the set Hermitian forms $L$ in
  Theorems~\ref{th:boundaryR} and Theorem~\ref{th:boundaryL} coincide.
  By Theorem~\ref{th:boundaryR}, $\sqminus L = \kappa$ if and only if
  there is a function $S(z)$ in $\SkappaCplus$ such that $b(x) =
  S(x)c(x)$ a.e.\ on~$\Delta$, that is, $f_0(x) - i = S(x)\big[ f_0(x)
  + i \big]$ or
  \begin{equation}
    \notag
     f_0(x) = i\, \frac{1 + S(x)}{1 - S(x)}
  \end{equation}
  a.e.\ on~$\Delta$.  The result thus follows from the correspondence
  \eqref{linearfrac} between $\SkappaCplus$ and $\Nkappa$.
\end{proof}

The preceding half-plane results can be recast in a different form
using Hilbert transforms.  The connection comes from equation
\eqref{keylimit2} and the formula \cite[p.~113]{RR1994}
\begin{equation}\label{projectionformula}
  Q_- = (I + iH)/2
\end{equation}
for the projection onto the orthogonal complement of $H^2$ in
$L^2(-\infty,\infty)$.  Here $H$ denotes the Hilbert transform, which
is defined for any $\varphi$ in $L^2(-\infty,\infty)$ by the
principal-value integral
\begin{equation}
  \notag
  (H\varphi)(x) = PV\, \frac{1}{\pi} \int_{-\infty}^\infty
  \frac{\varphi(t)}{t-x}\; dt 
  = \lim_{\epsilon \downarrow 0} \frac{1}{\pi} 
  \int_{|t-x|>\epsilon}   \frac{\varphi(t)}{t-x}\; dt .
\end{equation}
The limit exists pointwise a.e.\ on the real line and also in the
metric of $L^2(-\infty,\infty)$.  The compression of $H$ to
$L^2(\Delta)$ for any measurable subset $\Delta$ of $(-\infty,\infty)$
is denoted $H_\Delta$: if $\varphi \in L^2(\Delta)$, then
\begin{equation}
  \notag
  (H_\Delta \varphi)(x) = PV\; \frac{1}{\pi} \int_\Delta
  \frac{\varphi(t)}{t-x}\; dt
   = \lim_{\epsilon \downarrow 0} \, \frac{1}{\pi}
  \int_{|t-x|>\epsilon}   \frac{\varphi(t)}{t-x}\; dt
\qquad
\text{a.e.\ on} \; \Delta,
\end{equation}
where the limit is in the metric of $L^2(\Delta)$.  The operators $iH$
and $iH_\Delta$ are selfadjoint (in particular, $H_\Delta^* = -
H_\Delta$), and $iH$ is also unitary.  By the convolution theorem
(\cite[p.~35]{RR1985}, \cite[p.~111]{RR1994}), for any $\varphi,\psi
\in\, L^2(\Delta) \cap L^\infty(\Delta)$,
\begin{equation} \notag 
   H_\Delta \big[ \varphi (H_\Delta \psi) + (H_\Delta \varphi) \psi
  \big] = (H_\Delta \varphi)(H_\Delta \psi) - \varphi\psi
  \qquad \text{a.e.\ on}\;\; \Delta .
\end{equation}

\begin{lemma}\label{lemma}
    Let $f_0$ be a complex-valued measurable function defined on a Borel
  subset $\Delta$ of $(-\infty,\infty)$ such that $f_0 + i \neq 0$
  a.e.\ on $\Delta$.  Let $\calD$ be the set of measurable functions
  $\varphi$ on $\Delta$ such that $\varphi, f_0 \varphi \in
  L^2(\Delta)$.  If $\kappa$ is a nonnegative integer, there exists a
  function $f\in\Nkappa$ such that $f=f_0$ a.e.\ on $\Delta$ if and
  only if the Hermitian form
\begin{equation}\label{eq:bcform4}
L(\varphi,\psi) = \pi\, 
            {\ip{\big(H_\Delta - iI\big)(f_0\varphi)}{\psi}}_{L^2(\Delta)}
+ \pi\, {\ip{\varphi}{\big(H_\Delta - iI\big)(f_0\psi)}}_{L^2(\Delta)},
\end{equation}
$\varphi,\psi\in\calD$, has $\kappa$ negative squares.  The conclusion
remains true if $\calD$ is replaced by the set $\calE = \calD \cap
L^\infty(\Delta)$.
\end{lemma}

\begin{proof}
  Combining \eqref{keylimit2} and \eqref{projectionformula}, we obtain
  \begin{equation} \label{keylimit3}
    \lim_{\varepsilon\downarrow 0} \int_\Delta \int_\Delta
  \frac{\varphi(s)\overline{\psi(t)}}{s-t+i\varepsilon} \; ds\,dt
  = \pi\, {\ip{\big(H_\Delta - iI\big)\varphi}{\psi}}_{L^2(\Delta)},
   \qquad \varphi,\psi \in L^2(\Delta).
  \end{equation}
  It follows that the form $L$ in Theorem~\ref{th:boundaryL} is given
  by \eqref{eq:bcform4}, and so the first statement of the lemma
  follows from Theorem~\ref{th:boundaryL}.
  
  To prove the last statement, it is enough to show that for any
  $\varphi$ and $\psi$ in $\calD$, there are sequences $\varphi_n$ and
  $\psi_n$ in $\calE$ such that $L(\varphi_n,\psi_n) \to
  L(\varphi,\psi)$, since this implies that the number of negative
  squares of $L$ does not decrease when it is restricted to $\calE
  \times \calE$.  In fact, let $\varphi_n(t) = \varphi(t)$ or $0$
  according as $|\varphi(t)| \le n$ or $|\varphi(t)| > n$, and define
  $\psi_n(t)$ similarly ($n=1,2,\dots$).  Then $\varphi_n,\psi_n \in
  \calE$ for all $n$, and $\varphi_n\to\varphi$, $\psi_n \to \psi$ in
  the metric of $L^2(\Delta)$.  Hence also $f_0\varphi_n\to
  f_0\varphi$, $f_0\psi_n\to f_0\psi$ in the metric of $L^2(\Delta)$,
  and so $L(\varphi_n,\psi_n) \to L(\varphi,\psi)$ because $H_\Delta$
  is a continuous operator on $L^2(\Delta)$.
\end{proof}

We can convert the limit in \eqref{eq:bcform3} to a different form when
$f_0$ is real valued.

\begin{theorem}\label{th:realLoewner}
  Let $f_0$ be a real-valued measurable function defined on a Borel
  subset $\Delta$ of $(-\infty,\infty)$.  Let $\calD$ be the set of
  measurable functions $\varphi$ on $\Delta$ such that $\varphi, f_0
  \varphi \in L^2(\Delta)$.  If $\kappa$ is a nonnegative integer,
  there exists a function $f\in\Nkappa$ such that $f=f_0$ a.e.\ on
  $\Delta$ if and only if the Hermitian form
\begin{equation}\label{eq:bcform5}
L(\varphi,\psi) = 
\lim_{\varepsilon\downarrow 0} 
  \iint_{|t-s| > \epsilon}
  \frac{f_0(s) - f_0(t)}{s-t}
  \,\varphi(s)\overline{\psi(t)}\; ds\,dt,
\qquad \varphi,\psi\in\calD,
\end{equation}
has $\kappa$ negative squares
\end{theorem}

Integration in \eqref{eq:bcform5} is over the set
$(\Delta\times\Delta)\cap \{ (s,t) \colon |t-s| > \epsilon \}$.

\begin{proof}
  As in \cite[p.~37]{RR1985}, for any $\varphi,\psi$ in~$\calD$ we have
  \begin{align*}
    \lim_{\varepsilon\downarrow 0} 
  \iint_{|t-s| > \epsilon}
 & \frac{f_0(s) - f_0(t)}{s-t}
  \,\varphi(s)\overline{\psi(t)}\; ds\,dt \\[4pt]
  &= 
   \lim_{\varepsilon\downarrow 0} 
   \int_\Delta \bigg( \int_{|t-s| > \epsilon}
     \frac{f_0(s)\varphi(s)}{s-t} 
    \; ds\bigg) \overline{\psi(t)}  \,dt \\
&\hskip1cm   -  \lim_{\varepsilon\downarrow 0} 
    \int_\Delta \bigg( \int_{|t-s| > \epsilon}
     \frac{\varphi(s)}{s-t} 
    \; ds\bigg) f_0(t)\overline{\psi(t)}  \,dt \\[4pt]
  &= \pi\, {\ip{H_\Delta (f_0\varphi)}{\psi}}_{L^2(\Delta)}
       - \pi\, {\ip{H_\Delta \varphi}{f_0 \psi}}_{L^2(\Delta)}  
  \\[4pt]
  &= \pi\, {\ip{H_\Delta (f_0\varphi)}{\psi}}_{L^2(\Delta)} 
       + \pi\, {\ip{\varphi}{H_\Delta (f_0 \psi)}}_{L^2(\Delta)} .
  \end{align*}
  Since $f_0$ is real-valued, we can rewrite this in the form
  \eqref{eq:bcform4}, and therefore the assertion follows from
  Lemma~\ref{lemma}.
\end{proof}

In the definite case, Theorem~\ref{th:realLoewner} has a dual result
for functions having pure imaginary values \cite[pp.~39--41]{RR1985}.
A full generalization of the dual result to the class $\Nkappa$ is
unknown, but we are able to present some partial information.  For
simplicity, we assume that the given function $g_0$ and set $\Delta$
are bounded.

\begin{theorem}\label{th:dualresult}
  Let $g_0$ be a bounded real-valued measurable function on a bounded
  Borel subset~$\Delta$ of $(-\infty,\infty)$.  Let $\kappa$ be a
  nonnegative integer.
  \begin{enumerate}
  \item[(1)] There is a function $f\in\Nkappa$ such that $f= -
    H_\Delta g_0$ a.e.\ on $\Delta$ if and only if the Hermitian form
  \begin{equation}\label{eq:bcform6}
  L(\varphi,\psi) = \pi \int_\Delta \big[ \varphi \bar \psi
   -  (H_\Delta \varphi) (\overline{H_\Delta \psi}) \big] \; g_0 \; dt,
   \qquad \varphi,\psi\in  L^\infty(\Delta),
  \end{equation}
  has $\kappa$ negative squares.
\item[(2)] If the condition in $(1)$ holds and $g_0 \ge 0$ a.e.\ on
  $\Delta$, there is a function $h \in \Nkappaprime$, $\kappa' \le
  \kappa$, such that $h = ig_0$ a.e.\ on $\Delta$.
  \end{enumerate}
\end{theorem}

\begin{proof}
  (1) We apply Lemma~\ref{lemma} with the real-valued function $f_0 =
  -H_\Delta g_0$.  Since $g_0$ and $\Delta$ are bounded, $f_0 \in
  L^2(\Delta)$ and the set $\calE$ in Lemma~\ref{lemma} is $\calE =
  L^\infty(\Delta)$.  As in \cite[p.~39]{RR1985}, for any $\varphi$
  and $\psi$ in $L^\infty(\Delta)$, \eqref{eq:bcform4} is given by
\allowdisplaybreaks {
  \begin{align}
    L(\varphi,\psi) &= \pi\, 
            {\ip{\big(H_\Delta - iI\big)(f_0\varphi)}{\psi}}_{L^2(\Delta)}
+ \pi\, {\ip{\varphi}{\big(H_\Delta - iI\big)(f_0\psi)}}_{L^2(\Delta)}
          \notag  \\
&= \pi\, 
            {\ip{H_\Delta (f_0\varphi)}{\psi}}_{L^2(\Delta)}
+ \pi\, {\ip{\varphi}{H_\Delta (f_0\psi)}}_{L^2(\Delta)}
          \notag  \\
&= - \pi\, {\ip{f_0 \varphi}{H_\Delta \psi }}_{L^2(\Delta)}
    - \pi\, {\ip{H_\Delta \varphi}{f_0 \psi}}_{L^2(\Delta)} \notag \\
&= - \pi\, \int_\Delta \big[f_0 \varphi(H_\Delta \bar\psi)
      + (H_\Delta\varphi) f_0 \bar\psi    \big]\; dt \notag \\
&= \pi\, \int_\Delta \big[
\varphi  (H_\Delta \bar\psi) + (H_\Delta\varphi) \bar\psi
\big] H_\Delta g_0 \; dt \notag \\
&= -  \pi\, \int_\Delta \Big\{ H_\Delta \big[
\varphi  (H_\Delta \bar\psi) + (H_\Delta\varphi) \bar\psi
\big] \Big\} g_0 \; dt \notag \\
&= \pi\, \int_\Delta  \big[
\varphi\bar\psi - (H_\Delta\varphi)(\overline{H_\Delta \psi)}
 \big] g_0 \; dt . \notag
\notag 
\end{align}
We have shown that \eqref{eq:bcform4} coincides with
\eqref{eq:bcform6} on $\calE\times\calE$, and therefore the result
follows from Lemma~\ref{lemma}.
} 

(2) Assume the condition in (1) holds and $g_0 \ge 0$ a.e.\ on
$\Delta$.  Let $f_0 = - H_\Delta g_0$ a.e.\ on $\Delta$.  Then by the
first part of the theorem there is a function $f \in \Nkappa$ such
that $f(x+i0) = f_0(x)$ a.e.\ on $\Delta$.  Define $h(z)$ and $k(z)$
in the upper half-plane by
\begin{equation}
  \notag
  h(z) = k(z) + f(z) = \frac{1}{\pi} \int_\Delta \frac{g_0(t)}{t-z} 
  \; dt + f(z).
\end{equation}
Since $g_0 \ge 0$ a.e.\ on $\Delta$, $k(z)$ has nonnegative imaginary
part on the upper half-plane.  Hence the kernel
\begin{equation}
  \notag
  \frac{h(z) - \overline{h(\zeta)}}{\vphantom{t^{t^t}}z-\bar\zeta}
  = \frac{k(z) - \overline{k(\zeta)}}{\vphantom{t^{t^t}}z-\bar\zeta}
      + \frac{f(z) - \overline{f(\zeta)}}{\vphantom{t^{t^t}}z-\bar\zeta}
\end{equation}
has $\kappa' \le \kappa$ negative squares, that is, $h(z)$ belongs to
$\Nkappaprime$.  The boundary function $k_0(x) = k(x+i0)$ of $k(z)$
can be computed from \cite[Theorem A, p.~35]{RR1985}:
\begin{equation}
  \notag
  k_0(x) = i\, \big[ g_0(x) - i (H_\Delta g_0)(x)\big]
  \qquad    \; \text{a.e.\ on} \;\; \Delta.
\end{equation}
Therefore $h_0(x) = h(x+i0)$ is given by
\begin{equation}
  \notag
  \qquad
  h_0(x) = k_0(x) + f_0(x) = i\, \big[ g_0(x) - i (H_\Delta g_0)(x)\big] 
       - (H_\Delta g_0)(x) = ig_0(x)
\end{equation}
a.e.\ on $\Delta$.
\end{proof}

\section{Proof of the main result}\label{S:proof}

We use the Ball and Helton almost-commutant lifting theorem~\cite{BH}.
This result has been extended and refined in \cite{AADM2,MR97e:47030}.
Our statement follows \cite[Theorem 1.1]{AADM2}, which is slightly
more precise than the original.  In our applications, it is sufficient
to take the underlying Kre{\u\i}n spaces to be Hilbert spaces.  A
Hilbert space isometry $W \in \eL(\eG)$ is called an {\bf isometric
dilation} of an operator $T\in \lh$ if $\eG = \eH \oplus \eK$ for some
Hilbert space $\eK$, and 
\begin{equation} \notag W = \begin{bmatrix} T
& 0 \\ *&* \end{bmatrix} 
\end{equation} 
relative to the decomposition $\eG = \eH \oplus \eK$.  If $H$ is a
selfadjoint operator on a Hilbert space, $\indminus H$ is the
dimension of the spectral subspace of $H$ corresponding to the
interval~$(-\infty,0)$.

\begin{theorem}\label{th:ballhelton}
  For each $j=1,2$, let $T_j \in \eL(\eH_j)$ be a contraction on the
  Hilbert space $\eH_j$, let $W_j \in \eL(\eG_j)$ be an isometric
  dilation of $T_j$ on a Hilbert space $\eG_j$, and let $P_j$ be the
  projection of $\eG_j$ onto $\eH_j$.  Let $C \in \eL(\eH_1,\eH_2)$ be
  an operator such that $C T_1 = T_2 C$.
  \begin{enumerate}
  \item[(1)] If $\,\indminus (1 - C^*C) = \kappa$, there exists a
    pair $(\eE,\tildeC)$ such that $\eE$ is a closed $W_1$-invariant
    subspace of $\eG_1$ of codimension equal to $\kappa$ and $\tildeC$
    is a contraction operator on $\eE$ into $\eG_2$ satisfying
    \begin{equation}
      \notag
      \tildeC W_1 \vert_{\eE} = W_2 \tildeC
      \qquad \text{and} \qquad 
      P_2 \tildeC = C P_1 \vert_{\eE}.
    \end{equation}
  \item[(2)] If there is a pair $(\eE,\tildeC)$ with the properties in
    $(1)$, then $\indminus (1 - C^*C) \le \kappa$.
  \end{enumerate}
\end{theorem}

In part (1), the subspace $\eE$ can be chosen to be any closed
$W_1$-invariant subspace of $\eG_1$ of codimension equal to $\kappa$
such that $C P_1 \vert_{\eE}$ is a contraction \cite{AADM2}.

\begin{proof}[Proof of the Main Theorem] 
  Let $S$ be multiplication by~$z$ on $H^2$.  By parts (i) and (ii) of
  Definition~\ref{D:admissible}, the set of functions $h_{x'}(z) =
  \sum_{j=0}^\infty (A^jc,x') z^j$, $x' \in \calD$, is a linear
  subspace $\eH_c$ of $H^2$ which is invariant under $S^*$.  By part
  (iii) of Definition~\ref{D:admissible}, the formula
  \begin{equation}
    \notag
    X_0 \colon \sum_{j=0}^\infty (A^jc,x') z^j
    \to \sum_{j=0}^\infty (A^jb,x') z^j
  \end{equation}
  defines a linear operator $X_0$ from $\eH_c$ into $H^2$.
  
  (1) Assume $\calF$ has $\kappa$ negative squares.  We show that this
  assumption implies that $X_0 = Y_0 + K$, where $Y_0$ is a
  contraction and $K$ is a finite-rank operator, and hence $X_0$ is
  bounded in the norm of $H^2$.  To do this, let $\eX$ be $\eH_c$ as a
  vector space equipped with the linear and symmetric inner product
  \begin{equation}
    \notag
    {\ip{h_{x'}}{h_{y'}}}_{\eX} = \calF(x',y'), \qquad x',y' \in \calD.
  \end{equation}
  Since $\calF$ has $\kappa$ negative squares, $\eX$ contains a
  $\kappa$-dimensional subspace $\eN$ which is the antispace of a
  Hilbert space, and no $(\kappa+1)$-dimensional subspace of $\eX$ has
  this property~\cite[Lemma 1.1.1$^\prime$]{ADRS}.  Writing
  $$\eN^\perp = \{ h \in \eX \colon {\ip{h}{k}}_{\eX} = 0 \; \text{for
    all} \; k \in\eN \},$$
  we claim that
  \begin{equation}
    \label{nnperp}
    \eX = \eN + \eN^\perp 
    \quad \text{and} \quad\eN \cap \eN^\perp = \{0\}.
  \end{equation}
  If $h \in \eN \cap \eN^\perp$, then ${\ip{h}{h}}_{\eX} = 0$; this
  implies $h=0$ because $h \in \eN$ and the $\eX$-inner product is
  definite on $\eN$.  To see that $\eX = \eN + \eN^\perp$, choose a
  complete orthogonal set $h_1,\dots,h_\kappa$ in $\eN$ (in the
  $\eX$-inner product).  Then any $h\in\eX$ can be written
  \begin{equation}
    \notag
    h = \sum_{j=1}^\kappa 
         \dfrac{{\ip{h}{h_j}}_{\eX}}{{\ip{h_j}{h_j}}_{\eX}}\; h_j
         +
         \bigg(
           h - \sum_{j=1}^\kappa 
         \dfrac{{\ip{h}{h_j}}_{\eX}}{{\ip{h_j}{h_j}}_{\eX}}\; h_j
         \bigg).
  \end{equation}
  The first term here belongs to $\eN$ and the second is in
  $\eN^\perp$, so $\eX = \eN + \eN^\perp$ and the claim is proved.
  Thus $\eN^\perp$ is a subspace of $\eX$ having a finite-dimensional
  complement in $\eX$.  For any $h \in \eN^\perp$, by the definitions
  of $\calF$ and the $\eX$-inner product,
\begin{equation}
    \label{contractionpart}
    0 \le {\ip{h}{h}}_{\eX} = {\ip{h}{h}}_{H^2}
          - {\ip{X_0 h}{X_0 h}}_{H^2}.
\end{equation}
  Hence $X_0 \vert_{\eN^\perp}$ is a contraction in the norm of $H^2$.
  We are now able to construct operators $Y_0$ and $K$ with the
  required properties.  In fact, using \eqref{contractionpart} we
  first construct a contraction operator $Y_0$ on $\eH_c$ such that
  $Y_0 \vert_{\eN^\perp} = X_0 \vert_{\eN^\perp}$.  Then set $K=X_0 -
  Y_0$.  By construction $Y_0$ is a contraction, and by \eqref{nnperp}
  the range of $K$ is the range of $K\vert_\eN$, which has finite
  dimension; thus $X_0 = Y_0 + K$ where $Y_0$ is a contraction and $K$
  has finite rank.  In particular, $X_0$ is bounded in the norm of
  $H^2$.  Let $X$ be the extension by continuity of $X_0$ to the
  closure of $\eH_c$ in $H^2$.
  
  We apply Theorem~\ref{th:ballhelton} with these choices:
  \begin{enumerate}
  \item[] $\eH_1 = H^2$ and $T_1 = S$;
  \item[] $\eH_2 = \overline{\eH}_c$ (closure in $H^2$) and $T_2 =
    E_2^*SE_2$, where $E_2\colon \eH_2 \to H^2$ is inclusion;
  \item[] $\eG_1 = \eG_2 = H^2$ and $W_1 = W_2 = S$;
  \item[] $C = X^*$ (as an operator from $H^2$ to $\eH_2$).
  \end{enumerate}
  Trivially, $W_1$ is an isometric dilation of $T_1$.  We show that
  $W_2$ is an isometric dilation of~$T_2$.  Since $S^* \eH_c \subseteq
  \eH_c$, $S^* \eH_2 \subseteq \eH_2$.  Hence the matrix decomposition
  of $S^*$ relative to the decomposition $H^2 = \eH_2 \oplus
  \eH_2^\perp$ has the form
  \begin{equation}
    \notag
    S^* = 
    \begin{bmatrix}
      E_2^*S^*E_2 & * \\ 0 & *
    \end{bmatrix},
  \end{equation}
  which is equivalent to
\begin{equation}
    \notag
    W_2 = S = 
    \begin{bmatrix}
      E_2^*SE_2 & 0 \\ * & *
    \end{bmatrix}
    = \begin{bmatrix}
      T_2 & 0 \\ * & *
    \end{bmatrix}.
\end{equation}
  It is not hard to see that $CT_1=T_2C$ and $\indminus (1 - CC^*) =
  \indminus (1- X^*X) = \kappa$; since always $\indminus (1 - C^*C) =
  \indminus (1 - CC^*)$, we have $\indminus (1 - C^*C) = \kappa$.  By
  Theorem~\ref{th:ballhelton}(1), there is a pair $(\eE,\tildeC)$ such
  that
\begin{enumerate}
\item[(i)] $\eE$ is a closed subspace of $H^2$ of codimension~$\kappa$
  which is invariant under $S$,
\item[(ii)] $\tildeC:\eE\rightarrow H^2$ is a contraction,
\item[(iii)] $P_{\eH_2}\tildeC = X^*\vert_\eE$,
\item[(iv)] $\tildeC S \vert_\eE = S \tildeC$.
\end{enumerate}
For any function $F(z)$ in $\Szero$, write $\widetilde F(z) =
\overline{F(\bar z)}$.  By (i), $\eE = \widetilde B H^2$, where
$B(z)$ is a Blaschke product of degree~$\kappa$.  For
any $\varphi \in \Szero$, let $M_\varphi$ be the operator
multiplication by $\varphi$ on $H^2$.  We show that \smash{$\tildeC
M_{\widetilde B}$} commutes with $S$.  In fact, for each \smash{$h\in
H^2$}, $\widetilde B h\in\eE$ and
\begin{equation}
    \notag
    \tildeC M_{\widetilde B} Sh =
    \tildeC S (\widetilde B h) \stackrel{(iv)}{=\vphantom{X^x}} 
      S \tildeC\widetilde
    B h =S \tildeC M_{\widetilde B} h.
\end{equation}
  Since $\tildeC M_{\widetilde B}$ is a contraction by (ii), $\tildeC
  M_{\widetilde B} = M_{\tilde f}$ for some $f \in \Szero$.  For
  any $x'\in\calD$ and $h\in H^2$,
  \begin{multline}
    \!\!\!{\ip{\sum_{j=0}^\infty (A^jc,x') \,z^j}{M_{\tilde f}h}}_{\!H^2} =
  {\ip{\sum_{j=0}^\infty (A^jc,x') \,z^j}{\tildeC
      M_{\widetilde B}h}}_{\!H^2}  \\
  = {\ip{\sum_{j=0}^\infty (A^jc,x') \,z^j}{P_{\eH_2}\tildeC
      M_{\widetilde B}h}}_{\!H^2}
  \stackrel{(iii)}{=\vphantom{X^x}} {\ip{\sum_{j=0}^\infty (A^jc,x')
      \,z^j}{X^*M_{\widetilde B}h}}_{\!H^2} \\
  = {\ip{\sum_{j=0}^\infty (A^jb,x') \,z^j}{\widetilde B h}}_{\!\!H^2} .
  \notag
  \end{multline}
Choosing $h=1$, we obtain
\begin{equation} \label{E:nud1}
  {\ip{\sum_{j=0}^\infty (A^jc,x') \,z^j}{\tilde f}}_{H^2}
=
{\ip{\sum_{j=0}^\infty (A^jb,x') \,z^j}{\widetilde B}}_{H^2},
\qquad x'\in\calD.
\end{equation}
Thus \eqref{eq:twosums} and hence \eqref{eq:nud2} hold.  This
completes the proof of part (1) of the Main Theorem.

(2) Assume that \eqref{eq:nud2} holds with $f$ in $\Szero$ and $B$ a
Blaschke product of degree~$\kappa$.  We show that this assumption
implies that the operator $X_0$ defined at the beginning of the proof
is bounded.  Define \smash{$\tilde f$} and \smash{$\widetilde B$} as
in part (1) above.  The assumption \eqref{eq:nud2} can be rewritten in
the form \eqref{E:nud1}.  Hence for all $x'\in \calD$,
\begin{multline} 
  {\ip{z\tilde f}{\sum_{j=0}^\infty (A^jc,x') \,z^j}}_{H^2} =
  {\ip{\tilde f}{S^* \bigg( \sum_{j=0}^\infty (A^jc,x') \,z^j
      \bigg)}}_{H^2}\\
 ={\ip{\tilde f}{\sum_{j=0}^\infty (A^jc,A'x') \,z^j}}_{H^2}
  \stackrel{(\ref{E:nud1})}{=} {\ip{\widetilde B}{\sum_{j=0}^\infty (A^jb,A'x')
      \,z^j}}_{H^2} \\
= {\ip{z\widetilde B}{\sum_{j=0}^\infty (A^jb,x')
      \,z^j}}_{H^2}.  \notag
\end{multline} 
Writing as before $h_{x'}(z) = \sum_{j=0}^\infty (A^jc,x') \,z^j$ for
all $x'\in\calD$, we can extend this relation to
\begin{equation} \label{eq:ipid}
  {\ip{\tilde f h}{h_{x'}}}_{H^2} = {\ip{\widetilde B h}{X_0
      h_{x'}}}_{H^2}
\end{equation} 
first for every polynomial $h$ and then for every $h$
in $H^2$.  Now the boundedness of $X_0$ follows from the fact that
$\calE = \widetilde B H^2$ has codimension $\kappa$ in $H^2$.  Let $X$
be the extension by continuity of $X_0$ to the closure $\eH_2$ of
$\eH_c$ in $H^2$.  Since
\begin{equation*}
  \calF(x',y') = {\ip{h_{x'}}{h_{y'}}}_{H^2} 
- {\ip{Xh_{x'}}{Xh_{y'}}}_{H^2},
\qquad x',y'\in\calD,
\end{equation*}
to show that \eqref{eq:Fform} has $\kappa'$ negative squares for some
$\kappa' \le \kappa$, it is sufficient to show that $\indminus (1 -
X^*X)$ is at most~$\kappa$.  Write
\begin{equation} 
\notag 
X = \begin{bmatrix} X_1 \\ X_2 \end{bmatrix}
\end{equation}
relative to the decomposition $H^2 = \eE \oplus \eE^\perp$.  By
\eqref{eq:ipid}, $X^*M_{\widetilde B} = P_{\eH_2}M_{\tilde f}$ is a
contraction.  Therefore $X_1^* = X^*\vert_\eE$ and $X_1$ are
contractions.  Thus
\begin{equation} \notag 1 - X^*X = (1 - X_1^*X_1)
  - X_2^*X_2 ,
\end{equation} 
where $1 - X_1^*X_1 \ge 0$ and $- X_2^*X_2$ has rank at most~$\kappa$.
Hence $\indminus (1 - X^*X) \le \kappa$ and part (2) of the Main
Theorem follows.
\end{proof}

We restate the Main Theorem in a way which characterizes all solutions
of Nudel$'$man's problem in terms of the solutions of the
corresponding lifting problem.  In the statement of this result, we
identify two pairs $(f_1,B_1)$ and $(f_2,B_2)$ such that $f_1 = \gamma
f_2$ and $B_1 = \gamma B_2$ for some constant $\gamma$, $|\gamma| =
1$, and we adopt notation as in the proof of the Main Theorem.  In
particular, $S$ is multiplication by~$z$ on $H^2$.  For any $\varphi
\in \Szero$, let $M_\varphi$ be multiplication by $\varphi$ on $H^2$.
We write $\tilde \varphi(z) = \sum_0^\infty \bar \varphi_n z^n$ if
$\varphi(z) = \sum_0^\infty \varphi_n z^n$.

\begin{maintheorem_alt}
  Let $(A,b,c)$ be given data, $\calD$ an admissible set.  Define
  $\calF$ on $\calD\times \calD$ by \eqref{eq:Fform}, and let $\kappa$
  be a nonnegative integer.  Then
  \begin{equation}
    \notag
    \sqminus \calF \le \kappa
  \end{equation}
  if and only if there is a pair $(f,B)$, where $f \in \Szero$ and $B$
  is a Blaschke product of degree~$\kappa$, such that $f(A)c = B(A)b$.
  The set of all such pairs $(f,B)$ is in one-to-one correspondence
  with the set of all pairs $(\calE, \tildeC)$ such that
  \begin{enumerate}
  \item[(i)] $\calE$ is a closed subspace of $H^2$ of codimension
    $\kappa$ which is invariant under~$S$;
  \item[(ii)] $\tildeC \colon \calE \to H^2$ is a contraction operator;
  \item[(iii)] if $P_\calE$ is the projection onto $\calE$, then for
    all $x' \in\calD$,
    \begin{equation*}
      \tildeC^* \colon \sum_{j=0}^\infty (A^jc,x') z^j
    \to P_\calE \Big\{ \sum_{j=0}^\infty (A^jb,x') z^j\Big\};
    \end{equation*}  
  \item[(iv)] $\tildeC S|_\calE = S \tildeC$.
  \end{enumerate}
  The correspondence is determined by the relations $\calE = \widetilde B
  H^2$ and $\tildeC M_{\widetilde B} = M_{\tilde f}$.
\end{maintheorem_alt}

\begin{proof}
  If $\sqminus \calF = \kappa' \le \kappa$, then by part (1) of the main
  theorem there is a pair $(f_1,B_1)$ where $f_1 \in \Szero$ and $B_1$
  is a Blaschke product of degree~$\kappa'$, such that $f_1(A)c =
  B_1(A)b$.  If $\kappa' = \kappa$, the pair $(f,B) = (f_1,B_1)$ has
  the required properties.  Otherwise write $\kappa = \kappa' + r$ and
  let $(f(z),B(z)) = (z^rf_1(z),z^rB_1(z))$.  Then $f \in \Szero$ and
  $B$ is a Blaschke product of degree~$\kappa$.  If 
  \begin{equation}
    \notag
    f_1(z) = \sum_{j=0}^\infty \alpha_j z^j
    \quad \text{and}\quad 
    B_1(z) = \sum_{j=0}^\infty \beta_j z^j,
  \end{equation}
  then the coefficients of $f(z)$ consist of $r$ zeros followed by
  $\alpha_0,\alpha_1,\dots$, and similarly for $B(z)$.  In order to
  prove that $f(A)c = B(A)b$, we must therefore show that for all $x'
  \in \calD$,
  \begin{equation}
    \notag
    \sum_{j=0}^\infty \alpha_j (A^{j+r}c,x') 
     = \sum_{j=0}^\infty \beta_j (A^{j+r}b,x'). 
  \end{equation}
  By the invariance of $\calD$ under $A'$, this identity is equivalent
  to 
  \begin{equation}
    \notag
    \sum_{j=0}^\infty \alpha_j (A^jc,(A')^rx') 
     = \sum_{j=0}^\infty \beta_j (A^jb,(A')^rx'),
  \end{equation}
  which holds because $f_1(A)c = B_1(A)b$.  This proves the necessity
  part of the first statement.  Sufficiency follows from part (2) of
  the Main Theorem.
  
  Given a pair $(f,B)$ where $f \in \Szero$ and $B$ is a Blaschke
  product of degree~$\kappa$ such that $f(A)c = B(A)b$, define
  $(\calE, \tildeC)$ by $\calE = \widetilde B H^2$ and $\tildeC M_{\widetilde
    B} = M_{\tilde f}$.  Obviously (i) and (ii) hold, and we easily
  check (iv).  To prove (iii), we use the identity $f(A)c = B(A)b$ to
  show that 
  \begin{equation}
    \notag
    {\ip{\sum_{j=0}^\infty (A^jc,x')z^j}{\tilde f h}}_{\!\!H^2}
    = {\ip{\sum_{j=0}^\infty (A^jb,x')z^j}{\widetilde B h}}_{\!\!H^2}
  \end{equation}
  first for monomials $h(z) = z^k$ and hence for all $h \in H^2$.
  Therefore for all $x' \in \calD$ and $h \in H^2$,
  \begin{equation}
    \notag
    {\ip{\sum_{j=0}^\infty (A^jc,x')z^j}{\widetilde C \widetilde B h}}_{\!\!H^2}
    = {\ip{\sum_{j=0}^\infty (A^jb,x')z^j}{\widetilde B h}}_{\!\!H^2} ,
  \end{equation}
  and since $\calE = \widetilde B H^2$, this implies (iii).
  
  Conversely, if $(\calE, \tildeC)$ is a pair that satisfies
  (i)--(iv), the proof of the Main Theorem shows that $\calE = \widetilde
  B H^2$ and \smash{$\tildeC M_{\widetilde B} = M_{\tilde f}$}, where $f \in
  \Szero$ and $B$ is a Blaschke product of degree~$\kappa$, such that
  $f(A)c = B(A)b$.  The correspondence is one-to-one because the
  relation $\calE = \widetilde B H^2$ determines $B$ up to a constant
  factor of modulus one.
\end{proof}

\providecommand{\bysame}{\leavevmode\hbox to3em{\hrulefill}\thinspace}
\providecommand{\MR}{\relax\ifhmode\unskip\space\fi MR }
\providecommand{\MRhref}[2]{%
  \href{http://www.ams.org/mathscinet-getitem?mr=#1}{#2}
}
\providecommand{\href}[2]{#2}


\end{document}